\documentclass{amsart}

\usepackage{tikz}
\usetikzlibrary{matrix}

\usetikzlibrary{datavisualization}
\usetikzlibrary{datavisualization.formats.functions} 
\usetikzlibrary{calc}

\usepackage[T1]{fontenc}
\usepackage[latin1]{inputenc}
\usepackage{amsmath}
\usepackage{amssymb}
\usepackage{amsthm}
\usepackage{dsfont}
\usepackage{color}
\usepackage{graphicx}
\usepackage{graphics}
\theoremstyle{plain}\newtheorem{theorem}{Theorem}
\theoremstyle{plain}\newtheorem{corollary}{Corollary}
\theoremstyle{definition}
\theoremstyle{plain}\newtheorem{defi}{Definition}
\theoremstyle{plain}\newtheorem{lemma}[defi]{Lemma}
\theoremstyle{plain}\newtheorem{proposition}[defi]{Proposition}
\theoremstyle{definition}

\newcommand{\N}{\mathbb{N}}
\newcommand{\R}{\mathbb{R}}

\begin{document}
\title[Homogenization of a diffusion driven by growing incompressible drift]{A generic approach to homogenization of a diffusion driven by growing incompressible drift}

\author[Brice Franke, Shuenn-Jyi Sheu]{Brice Franke, Shuenn-Jyi Sheu}

\address{Univ Brest, UMR CNRS 6205, Laboratoire de Math\'ematiques de Bretagne Atlantique, 6 Avenue Le Gorgeu, CS 93837, 29238 Brest, Cedex 3, France}
\email{brice.franke@univ-brest.fr}
\address{Department of Mathematics, National Central University, 
No.300, Jhongda Rd., Jhongli City, Taoyuan County 32001, Taiwan (R.O.C.)}
\email{sheusj@math.ncu.edu.tw}
\date{\today}

\begin{abstract} 
We study how the resolvent-family of a diffusion behaves, as the drift grows to infinity. 
The limit turns out to be a selfadjoint pseudo-resolvent. After reduction of the underlying Hilbert-space,
this pseudo-resolvent becomes a resolvent to a strongly continuous semi-group of contractions.
We prove that this semi-group is associated to some Hunt-process on some suitable state-space which is constructed from equivalence classes of the drifts trajectories. 
Finally, we show a distributional limit theorem for the accelerated diffusion toward the associated Hunt process.
\end{abstract}

\maketitle

\section{Introduction}

The influence of a large incompressible drift on a diffusion is a subject which has gained some attention in recent years.
One issue is to understand how fast the diffusion converges toward its equilibrium when the drift becomes larger.
This can be studied in terms of the operator-norm of the semigroup (see Constantin, Kislev, Ryzhik, Zlatos (2008)), in terms of the spectral gap of the 
generator (see Franke, Hwang, Pai, Sheu (2010)) or in terms of asymptotic variance (see Hwang, Normand, Wu (2015), Duncan, Leli\`evre, Pavliotis (2016), Franke, Hwang, Ouled-Said, Pai (2021)). \\
Another issue is to search, as the drift grows to infinity, for homogenization along the trajectories of the flow underlying the drift. 
In this situation the further and further growing drift leads to a diffusion which in the limit propagates along the flow lines at practically infinite speed, while the diffusive path length, which is responsible for moving from one flow line to another, stays bounded. 
In the process of growing drift the position on the flow trajectory of such a diffusion becomes highly volatile and the only sustainable information is the account for the successive flow lines which have been visited by the diffusion.
Therefore, increasing the speed of the drift to infinity should lead to a limit diffusion which takes its values on a suitable state space of trajectories.
However, this picture is too rough, since one has to account for the possibility that two trajectories come arbitrary close to each other on large time scales. 
In that situation an infinitesimal diffusive increment suffices to move from one trajectory to another. 
One thus has to use the dynamics of the underlying flow to introduce suitable equivalence classes on the underlying state space which account for this possibility.
In this manuscript it is our intention to construct a suitable state space for the limit diffusion.
It will turn out that the right state space is a quotient space with respect to a non-local equivalence relation and thus the limit diffusion might not be describable through some differential operator anymore. 
It is then natural to use a resolvent family or some Dirichlet form to describe those diffusions. 
Our investigation does not limit itself to two dimensions.\\
To our knowledge most of the existing results on averaging along the trajectories of a drift accelerated diffusion deal with the two dimensional situations (see Hwang and Sheu (2000), Freidlin and Weber (2004), Ishi and Souganidis (2012), Dolgopyat, Freidlin and Koralov (2012). 
A paper by Barret and von Renesse (2014) also deals with higher dimensional situations. 
In this paper they use advanced Dirichlet form theory to analyse the asymptotic behavior of the drift accelerated diffusion through some observables, which are first integrals of the accelerating drift. 
Their situation is more general, since they do not assume incompressibility of the accelerating flow with respect to some measure of stationarity for the diffusion. 
However, it has to be noted, that the construction of their limit process also uses information from the level sets of the observable and thus is not an generic process constructed only from the geometry of the driving flow and the underlying manifold.\\
In this article we want to present a generic construction for a quotient space and a sole diffusion on it, which underlies all the limit processes that one can obtain from observation of the accelerating diffusion through first integrals.
We think that this might help to understand the underlying mechanisms which drive homogenisation along trajectories.  \\
Our starting point is to understand how the resolvent of the diffusion behaves as the drift grows to infinity.
We will see that those resolvent operators converge strongly toward a self-adjoint pseudo-resolvent family (see Theorem \ref{L5} and Theorem \ref{PseudoR}). It will turn out that the range of this pseudo-resolvent is the closure of the kernel of the drift-generating first order differential operator. 
Once restricted to this sub-space, the pseudo-resolvent becomes a true resolvent family. 
In order to find a Dirichlet form and a Hunt process for this resolvent we have to further investigate the closure of the kernel associated with the drift vectorfield. 
It turns out that it contains a dense set of smooth drift-invariant functions (see Theorem \ref{regularity}). 
This will be essential to prove that the range of the pseudo resolvent is isometric to a Hilbert space of square integrable functions on a suitable quotient space (see Theorem \ref{Isometrie}). 
This enables us to give explicit representations for the Dirichlet form of the limit diffusion (see Theorem \ref{Dform}). 
Some standard arguments show that the semi groups of the accelerated diffusions converge toward the semi group asociated to this Dirichlet form (see Theorem \ref{HalbgruppenKonvergenz}).  Finally we also prove a weak limit theorem for observables of the associated stochastic diffusions toward the Hunt process which is associated to the above Dirichlet form (see Theorem \ref{C2Konvergenz}). \\

\section{The mathematical setting}

Let $ (M,g) $ be a $ d $-dimensional oriented compact Riemannian manifold without boundary. We denote the metric distance of two points 
$ x,y\in M $ by $ d_M(x,y) $. For the canonical volume element on $ M $ we will use the notations $ {\rm vol}(dx) $ (see Chavel (1984), p.6 for a definition). 
Let $ U $ be a function from $ C^1( M) $  such that $ \pi(dx)=\exp(-U(x)){\rm vol}(dx) $ turns out to be a probability measure on $  M $. 
By construction the measure $ \pi $ has full support on the manifold $ M $.
We will denote by $ L^2(M,\pi) $ the Hilbert space of square integrable real valued functions with inner product defined for $ f,h\in L^2(M,\pi) $ as 
$$   \langle f,h\rangle  :=\int_M fhd\pi .$$
We use a family of local coordinates  $ \Psi_\alpha :V_\alpha\rightarrow\R^d;\ \alpha\in I $ to introduce the Laplace operator 
$$  \Delta : C^2(M)\rightarrow C(M); f\mapsto \frac{1}{\sqrt{g}}\sum_{k,j}^d\partial_{x_j}\Big(g^{jk}\sqrt{g}\partial_{x_k}f\Big), $$
where as usual in Riemannian geometry $ \sqrt{g} $ denotes the square root of the determinant of the matrix $ G $ with entries 
$ g_{jk}=g(\partial_{x_j},\partial_{x_k}) $ and where $ g^{jk} $ denote the entries of the inverse matrix $ G^{-1} $
(see Chavel (1984), p.5). 
Note however, that with the above definition the Laplace operator is a differential operator with spectrum contained in $ ]-\infty,0] $.
We denote by $ (A,{\rm Dom}(A)) $ the Friedrichs extension in $ L^2(M,\pi) $ of the differential operator 
$$   A_o :C^\infty( M)\rightarrow C^\infty( M);f\mapsto\Delta f-\nabla U\cdot\nabla f  .$$
We assume that the operator $ A $ has a discrete spectrum $$ 0=\lambda_0>-\lambda_1\geq -\lambda_2\geq..., $$ with 
$ \lambda_n\rightarrow\infty $ as $ n\rightarrow\infty $ and that the corresponding eigenvectors $ \phi_0,\phi_1,\phi_2,... $ form a complete orthonormal basis of $ L^2(M,\pi) $. 
Note that $ \phi_0 $ spans the constant functions on $ M $.
Then every element $ \varphi\in L^2(M,\pi) $ can be represented as 
$$ \varphi=\sum_{k=0}^\infty a_k\phi_k  $$ with suitable coefficients $ a_k,k\in\N $.
We use the orthonormal basis $ (\phi_k)_{k\in\N} $ to introduce the Sobolev spaces 
$$ H^m :=\left\{\varphi=\sum_{k=1}^\infty a_k\phi_k;\ \sum_{k=1}^\infty|a_k|^2\lambda_k^m<\infty\right\}.$$
For notational reasons we will denote $ H^0 $ by $ H $ in the following. 
Once restricted to the orthogonal complement of the constant functions the domain of the unbounded operator $ A $ is $ H^2 $ and the domain of its associated Dirichlet form
$$ {\mathcal E}(f,h) :=\langle f,h\rangle_{H^1}=\int_M\nabla f\cdot\nabla hd\pi $$
in $ H^1 $. When $ f $ is in $ H^2 $, then one further has
$$  {\mathcal E}(f,h)=\langle f,h\rangle_{H^1}=-\int_M hAfd\pi  .$$
In the previous expression we used the gradient of a function $ f: C^1(M)\rightarrow\mathbb{R} $, which is defined in local coordinates as a tangential vector with $ j $-th component
$$ (\nabla f)_j=\sum_{i=1^d}g^{ij}(\partial_{x_i}f)\partial_{x_j} .$$
Here and in the following, we use the lightened notation $ v\cdot w:=g(v,w) $ for tangent vectors $ v,w \in T_xM $. 
Let $ b $ be a $ C^1 $-vectorfield on $  M $ with the properties that $ {\rm div}(e^{-U}b)=0 $ and 
$ H^1\subset {\rm Dom}(b\cdot\nabla) $. In local coordinates the divergence of a $ C^1 $-vectorfield $ b=\sum_j^d\eta^j\partial_{x_j} $ is defined through
$$ {\rm div}(b) :=\frac{1}{\sqrt{g}}\sum_j^d\partial_{x_j}\Big(\eta^j\sqrt{g}\Big) $$
(see Chavel (1984) p.5). With those definitions and conventions we have
\begin{eqnarray*}
   \int_M f\Big(\Delta-\nabla U\nabla\Big) hd\pi = -\int_M\nabla f\cdot\nabla h d\pi   \ \ \ \ \forall  f,h\in C^2(M) 
\end{eqnarray*}
and
\begin{eqnarray*}
     \int_Mfb\cdot hd\pi=-\int_Mhb\cdot\nabla fd\pi   \ \ \ \ \forall  f,h\in C^1(M).
\end{eqnarray*}
For $ c\in\R $ we then define the differential operators
$$ A_{c,o} :C^\infty( M)\rightarrow C^\infty( M);f\mapsto Af+cb\cdot\nabla f .$$

\begin{proposition}
There exists a closure $ (A_c,{\rm Dom}(A_c)) $ of $ (A_{c,o},C^\infty( M)) $ such that that 
$ {\rm Dom}(A_c)=H^2={\rm Dom}(A) $ and for all $ \psi\in{\rm Dom}(A_c) $
$$   A_c\psi=A\psi+cB\psi ,$$
where $ B $ is an antisymmetric operator, which equals $ f\mapsto b\cdot\nabla f $ when restricted to 
$ C^\infty( M) $.\footnote{It might be better to restate Proposition 1 as follows : The operator $ (A_c,H^2) $ is a closure of $ (A_{c,o},C^\infty( M)) $ such that that for all $ \psi\in H^2 $
$$   A_c\psi=A\psi+cB\psi ,$$
where $ B $ is an antisymmetric operator, which equals $ f\mapsto b\cdot\nabla f $ when restricted to 
$ C^\infty( M) $.}
\end{proposition}
{\bf Proof :}
Since for constant functions the statement is obvious, we restrict our consideration to functions $ \psi $ from the orthogonal complement of the constant functions in $ L^2(M,\pi) $.
For a sequence $ \psi_n,n\geq1 $ of functions from $ C^\infty(M) $, which converges in $ H^1 $ toward an element $ \psi\in H^1 $, 
we know that 
$$ \int_M|b\cdot\nabla(\psi_n-\psi_m)|^2d\pi\leq\sup_{x\in M}g(b,b)(x)\|\psi_n-\psi_m\|_{H^1}^2\longrightarrow0 .$$
Therefore, $ b\cdot\nabla\psi_n,n\geq1 $ is a Cauchy sequence in $ H $. Therefore, the operator
$ B\psi :=\lim b\cdot\nabla\psi_n $ is well defined for all $ \psi\in H^1 $.
Moreover, if $ \varphi $ is an element from $ H^1 $, then there exists a sequence $ \varphi_n,n\geq1 $ 
in $ C^\infty( M) $ such that $ \varphi_n\rightarrow\varphi $ in $ H^1 $. It then follows that
$$ \langle\varphi,B\psi\rangle=\lim_{n\rightarrow\infty}\int_M\varphi_nb\cdot\nabla\psi_nd\pi    
         =-\lim_{n\rightarrow\infty}\int_M\psi_nb\cdot\nabla\varphi_nd\pi=-\langle\psi,B\varphi\rangle .$$
This implies that the operator $ (B,H^1) $ is anti-symmetric.\\[1mm]
We already mentioned that $ {\rm Dom}(A)=H^2 $.
For $ \varphi\in H^2 $ we define
$$  A_c\varphi :=A\varphi+cB\varphi .$$ 
We have to show that $ (A_c,H^2) $ is a closed extension of the differential-operator 
$ (A_{c,o},C^\infty( M)) $.
Let $ \psi_n,n\geq0 $ be a sequence from $ H^2 $ and $ \psi,g\in H $ such that one has 
$$ \psi_n\longrightarrow\psi \ \ {\rm and} \ \  A_c\psi_n\longrightarrow g \ \ {\rm in}\  H .$$ 
We need to prove that $ \psi\in H^2 $ and $ A_c\psi=g $.  
Note that 
$$  \|\psi_n-\psi_m\|_{H^1}=-\langle \psi_n-\psi_m,A_c(\psi_n-\psi_m)\rangle\longrightarrow0 \ \ 
  {\rm as}\ n,m\rightarrow\infty .$$
This implies that $ \psi_n,n\geq1 $ is a Cauchy-sequence in $ H^1 $. 
Therefore we have that $ B\psi_n,n\geq1 $  converges in $ H $ toward an element $ h\in H $ and we have that $ B\psi=h $. 
Therefore, 
$$  A\psi_n=A_c\psi_n-cB\psi_n\longrightarrow g-ch \ \ {\rm as}\ n\rightarrow\infty.$$
Since $ A $ is a closed operator, we obtain that $ \psi $ is in $ H^2 $ and $ A\psi=g-ch $.
This then implies that $ \psi\in{\rm Dom}(A)=H^2 $ and $ A_c\psi=A\psi+cB\psi=g $.
\hfill $ \Box $\\[1ex]
The next proposition is a version of Rellich's lemma for our setting.

\begin{proposition}\label{Rellich}
If a sequence of elements $ \psi_n,n\geq1 $ from $ H^1 $ is bounded with respect to the $ H^1 $-norm, then there
exists an element $ \psi\in H^1 $ and a sub-sequence  $\psi_{n_k}; k\geq 1 $ which converges toward $ \psi $ weakly in $ H^1 $ and strongly in $ H $. 
\end{proposition}
{\bf Proof :} See Appendix for a proof.
\hfill $ \Box $

\section{The convergence of the resolvent}

Let $ \psi_c\in H^1 $ be the unique soultion of the equation
\begin{eqnarray}   \label{PsiEquation}
A_c\psi-\lambda\psi=g . 
\end{eqnarray}
The resolvent $ R_\lambda^{(c)} $ associated to the operator $ (A_c,{\rm Dom}(A_c)) $ is the family of bounded operators defined through
\footnote{Our definition corresponds to the definition stated on p.209 from Yoshida, Functional Analysis, Springer Grundlehren der mathematischen Wissenschaften, Springer Verlag Berlin and also on p.15 of Fukushima, Dirichlet Forms and Markov Processes, North Holland, Amsterdam} 
$$   R_\lambda^{(c)}g=-\psi_c;\ \lambda>0 .$$
This can be restated as
$$  R_\lambda^{(c)}g=-(A_c-\lambda)^{-1}g .$$
We want to analyze the limit behavior of the family of resolvents $ (R_\lambda^{(c)})_{\lambda\geq0} $ as $ |c|\rightarrow\infty $. This will later help us understanding the asymptotic behavior of the semigroups $ (T_t^{(c)})_{t\geq0} $ that are generated by $ (A_c,{\rm Dom}(A_c)) $ as $ |c|\rightarrow\infty $.
We will use ideas that where developed by Bhattacharya, Gupta and Walker (1989) for the asymptotic analysis of limiting diffusion coefficients in periodic homogenization. In the present context, we make use of their method to obtain information on the limit behavior of the resolvent family $ R_\lambda^{(c)} $ as $ |c|\rightarrow\infty $. We will prove the following result.

\begin{theorem}\label{L5}
For $ g\in H $, $ \lambda>0 $ and as $ |c|\rightarrow\infty $ the functions $ R_\lambda^{(c)}g $
converge in $ H^1 $ toward a limit of the form 
$$    R^\ast_\lambda g:=-P_{1,\lambda}(A-\lambda)^{-1}g , $$
where $ P_{1,\lambda} $ is the orthogonal projection onto the space 
$$ K:={\rm Ker}(B)\cap H^1=\Big\{\varphi\in H^1: B\varphi=0\Big\} $$ 
with respect to the scalar product which is defined on $ H^1 $ through
$$ \langle\psi,\varphi\rangle_{1,\lambda} :=\langle\psi,\varphi\rangle_{H^1}+\lambda\langle\psi,\varphi\rangle  .$$
\end{theorem}

{\bf Remark :} Later in Proposition \ref{OrthoKern} we will see that $ R_\lambda^\ast g=0 $ when $ g\in{\rm Range}(B)\cap H^1 $.\\[0.5mm]

{\bf Proof  of Theorem \ref{L5} :}
If we denote $ \psi_c :=-R_\lambda^{(c)}g $, then one has from (\ref{PsiEquation}) that
\begin{eqnarray*}\label{RsGl1}
   A\psi_c+cB\psi_c-\lambda\psi_c=g .
\end{eqnarray*}
Note that the operator $ A-\lambda $ can be inverted and this can be restated as
\begin{eqnarray}\label{RsGl2}
  (1+c(A-\lambda)^{-1}B)\psi_c=(A-\lambda)^{-1} g= :\gamma.
\end{eqnarray}
It was noted by Bhattacharya, Gupta and Walker (1989), that the operator 
$$ L :=(A-\lambda)^{-1}B : \ H^1\rightarrow H^2 \hookrightarrow H^1 $$ is an antisymmetric compact operator on the space $ H^1 $ with scalar product
$$  \langle\psi,\varphi\rangle_{1,\lambda} :=\langle\psi,\varphi\rangle_{H^1}+\lambda\langle\psi,\varphi\rangle 
=\langle\psi,(\lambda-A)\varphi\rangle .$$
As a consequence there exist eigenvalues $ (i\mu_k)_{k\in\N} $ and
corresponding eigenfunctions $ (\theta_k)_{k\in\N} $ with the properties  :\\[0.5mm]
i) One has $ \mu_k\in\R\backslash\{0\} $ for all $ k\in\N $ and 
$ \lim_{k\rightarrow\infty}\mu_k=0 $.\\[0.5mm]
ii) The family $ (\theta_k)_{k\in\N} $ is a complete orthonormal basis (ONB) for $ H^1\cap K^\perp $, where 
\begin{eqnarray} \label{KSpace}
 K=\Big\{\varphi\in H^1;Lh=0\Big\}={\rm Ker}\Big((A-\lambda)^{-1}B\Big)={\rm Ker}(B)\cap H^1 .
 \end{eqnarray}
iii) For all $ f\in H^1 $ one has 
$$     f=f_K+\sum_{k=1}^\infty\alpha_k\theta_k ,$$
where $ f_K\in K $ is the orthogonal projection of $ f $ toward $ K $ with respect to the scalar
product $ \langle\cdot,\cdot\rangle_{1,\lambda} $ and $ \alpha_k:=\langle f,\theta_k\rangle_{1,\lambda} $.\\[0.5mm]
If we represent $ \psi_c $ and $ \gamma :=(A-\lambda)^{-1}g $ in the ONB mentioned above, we obtain
$$   \psi_c=\psi_{c,K}+\sum_{k=1}^\infty\alpha_{c,k}\theta_k \ \ \ {\rm and} \ \ \ 
                        \gamma=\gamma_K+\sum_{k=1}^\infty\beta_k\theta_k .$$
The equation (\ref{RsGl2}) then implies the following relations
$$   \psi_{c,K}=\gamma_K \ \ \ {\rm and} \ \ \ \alpha_{c,k}=\frac{\beta_k}{1+ic\mu_k} .$$
Therefore, as $ |c|\rightarrow\infty $
$$  \psi_c=\gamma_K+\sum_{k=1}^\infty\frac{\beta_k}{1+ic\mu_k}\theta_k
                                 \stackrel{H^1}{\longrightarrow}\gamma_K= :\psi_\ast .$$
This finishes the proof. 
\hfill $ \Box $

\vspace{1ex}

As a consequence of Theorem \ref{L5} we have that for every $ \lambda>0 $ and $ g\in H $ the functions $ R_\lambda^{(c)}g $ converge in $ H^1 $ toward
\begin{eqnarray}  \label{ResolventenAusdruck}  
R_\lambda^\ast g=- \psi_\ast=-P_{1,\lambda}(A-\lambda)^{-1}g \in K \ \ \ \mbox{as $ c\rightarrow\infty $}.
\end{eqnarray}
In the following section, we will see that the famly of linear operators $ R_\lambda^\ast;\lambda>0 $ form a pseudo-resolvent family on $ H $. \\[0.5mm]

{\bf Remark :} If we define  $ \gamma=(A-\lambda)^{-1}g\in H^1 $, then we know from Theorem \ref{L5} that
$$ \psi_\ast=P_{1,\lambda}(A-\lambda)^{-1}g=P_{1,\lambda}\gamma .$$ 
This means that $ \psi_\ast $ is the unique minimizer in $ {\rm Ker}(L)={\rm Ker}(B)\cap H^1 $ of the functional
$$ \psi\mapsto \|\psi-\gamma\|_{1,\lambda}  .$$
Thus, we obtain the decomposition
$ \psi_\ast=\varphi_\ast+\gamma $, where $ \gamma $ is the unique solution of 
$$ (A-\lambda)\gamma=g $$
and $ \varphi_\ast:=\psi_\ast-\gamma $ is the unique minimizer of the functional 
$$  \|\varphi\|_{1,\lambda}^2
              =\|\varphi\|_{H^1}^2+\lambda\|\varphi\|^2={\mathcal E}_\lambda(\varphi)  $$ 
restricted to the set of $ \varphi\in H^1 $ satisfying $ \varphi+\gamma\in{\rm Ker}(B) $, which is equivalent to 
$$ B\varphi=-B\gamma .$$

\vspace{1ex}

For later use we now generalize some of the concenpts from the proof of Theorem \ref{L5}. 
We first observe that for all $ m\in\N\backslash\{0\} $ the operators $ (A-\lambda)^{-m}B $ 
are antisymmetric compact operators on the Hilbert space $ H^m $ with scalar product 
$$  \langle \psi,\varphi\rangle_{m,\lambda} :=\langle \psi,(\lambda-A)^m\varphi\rangle . $$
For two subspaces $ V, W $ in $ H^m  $ we note $ V\perp_{m,\lambda}W $ if $ V $ is orthogonal to $ W $ with respect to the scalar product $ \langle.,.\rangle_{m,\lambda} $.\\[0.5mm]
For a given $ g\in H $, we say that a function $ \varphi \in H $ is a weak solution for the equation 
$$  B\varphi\stackrel{w}{=}g $$
iff for all $ \psi\in H^1 $ holds
$$  \langle \varphi ,B\psi\rangle=-\langle g,\psi\rangle  .$$

In the following, we denote by $ \overline{{\rm Ker}(B)\cap H^m} $ the closure of the set $ {\rm Ker}(B)\cap H^m $ in $ H $.

\begin{theorem} \label{regularity}
For all $ m\in\N\backslash\{0\} $ the closure of the set $ {\rm Ker}(B)\cap H^m $ in $ H $ is just the set of weak solutions of the equation $ B\psi=0 $; i.e. :
$$ \overline{{\rm Ker}(B)\cap H^m}=\Big\{\varphi\in H :\langle\varphi,B\psi\rangle=0, \ \forall \psi\in H^1\Big\}   .$$
\end{theorem}
{\bf Proof :}
First we note that for an element $ \varphi\in\overline{{\rm Ker}(B)\cap H^m} $ there exists some sequence of elements 
$ \varphi_n\in H^m $ with $ \varphi_n\rightarrow\varphi $ in $ H $ and $ B\varphi_n=0 $ for all $ n\in\N $. Therefore, for all 
$ \psi\in H^1 $ we have that
$$   \langle\varphi,B\psi\rangle=\lim_{n\rightarrow\infty}\langle\varphi_n,B\psi\rangle
                                =-\lim_{n\rightarrow\infty}\langle B\varphi_n,\psi\rangle = 0 .$$
This proves that $ \varphi $ is a weak solution of the equation $ B\varphi=0 $.
Thus it remains to prove that all weak solutions $ \varphi\in H $ of $ B\varphi=0 $ are in the closure of the set 
$ {\rm Ker}(B)\cap H^m $. If we assume that there exists a $ \tilde{\varphi}\in H $ with $ B\tilde{\varphi}\stackrel{w}{=}0 $ which is 
not in the closure of the set $ {\rm Ker}(B)\cap H^m $ then there exists a $ \psi\in H $ which satisfies
$ \langle \tilde{\varphi},\psi\rangle>0 $ and
\begin{eqnarray} \label{Kernrelation}
 \langle \varphi,\psi\rangle=0 \ \ \mbox{for all} \ \varphi\in{\rm Ker}(B)\cap H^m={\rm Ker}((A-1)^{-m}B)\cap H^m .
\end{eqnarray}
Note that since the operator $ B $ maps $ H^{m} $ to $ H^{m-1} $ we have that the operator
$$  (A-1)^{-m}B :H^m\rightarrow H^{3m-1}\hookrightarrow H^m $$
is a compact antisymmetric operator on $ H^m $ with respect to the the scalar product $ \langle.,.\rangle_{m,\lambda} $. 
The relation in (\ref{Kernrelation}) implies 
$$ \left\langle\varphi,(A-1)^{-m}\psi\right\rangle_{m,1}=\left\langle\varphi,(1-A)^m(A-1)^{-m}\psi\right\rangle=(-1)^m\langle\varphi,\psi\rangle=0 $$
for all $ \varphi\in{\rm Ker}((A-1)^{-m}B)\cap H^m $; i.e. :
$$  (A-1)^{-m}\psi\perp_{m,1}{\rm Ker}((A-1)^{-m}B) .$$
This implies by the closed range theorem (see \cite{Yoshida})
$$  (A-1)^{-m}\psi\in{\rm Range}((A-1)^{-m}B) ,$$
and therefore there exists a $ g\in H^m $ with the property 
$$  (A-1)^{-m}\psi=(A-1)^{-m}Bg $$
This means that $ \psi=Bg $.
The fact that $ \tilde{\varphi} $ satisfies $ B\tilde{\varphi}\stackrel{w}{=}0 $ then implies that
$$ \langle \tilde{\varphi},\psi\rangle=\langle \tilde{\varphi},Bg\rangle=0 .$$
This is a contradiction to our starting assumption. Thus $ \tilde{\varphi} $ must be in the closure of 
$ {\rm Ker}(B)\cap H^m $.
\hfill $ \Box $

\begin{corollary} \label{Dicht}
For all $ k\in\N $ there exists in the set of weak solutions of $ B\varphi=0 $ a dense set of $ C^k $-functions; i.e. :
$ \left\{\varphi\in H :\langle\varphi,B\psi\rangle=0, \ \forall \psi\in H^1\right\}\cap C^k(M) $ is dense in
$ \left\{\varphi\in H :\langle\varphi,B\psi\rangle=0, \ \forall \psi\in H^1\right\} $. 
\end{corollary}
{\bf Proof :}
This follows from the previous proposition and the Sobolev embedding theorems, which states that $ H^m\subset C^k(M) $ for $ m $ sufficiently large.
\hfill $ \Box $
\\[1ex]
For later use we also need some generalization of Theorem \ref{regularity} to higher order Sobolev spaces.
On the space $ H^1 $ we introduce the norm 
$$    \|\psi\|_{1,\lambda} :=\sqrt{\langle\psi,\psi,\rangle_{1,\lambda}}     .$$
For $ m\in\N\backslash\{0\} $ the closure of the set $ {\rm Ker}(B)\cap H^m $ with respect to the norm $ \|.\|_{1,\lambda} $ by $ \overline{{\rm Ker}(B)\cap H^m}^{(1,\lambda)} $.

\begin{proposition} \label{H1dicht}
For all $ m\in\N\backslash\{0\} $ the set $ \overline{{\rm Ker}(B)\cap H^m}^{(1,1)} $ is just the set of weak solutions in $ H^1 $ of the equation $ B\psi=0 $; i.e. :
 $$                  \overline{{\rm Ker}(B)\cap H^m}^{(1,1)} 
          =\left\{\varphi\in H^1 :\langle\varphi,B\psi\rangle=0, \forall\psi\in H^1\right\}={\rm Ker}(B)\cap H^1 .$$
\end{proposition}
{\bf Proof :}  For an element 
$ \varphi\in  \overline{{\rm Ker}(B)\cap H^m}^{(1,1)}  $ there exists some sequence of 
elements $ \varphi_n\in H^m$ satisfying $ \|\varphi_n-\varphi\|_{1,1}\rightarrow0 $ and $ B\varphi_n=0 $.
Therefore, for all $ \psi\in H^1 $ we have that
$$ \langle\varphi,B\psi\rangle = \lim_{n\rightarrow}\langle\varphi_n,B\psi\rangle
    =-\lim_{n\rightarrow\infty}\langle B\varphi_n,\psi\rangle=0 .$$
Again we have to prove that all weak solutions $ \varphi\in H^1 $ of $ B\varphi=0 $ are in the
closure of the set $ {\rm Ker}(B)\cap H^m $ with respect to the norm $ \|.\|_{1,1} $. Assuming the existence of a $ \tilde{\varphi}\in H^1 $ with $ B\tilde{\varphi}\stackrel{w}{=}0 $ which is not in the closure of $ {\rm Ker}(B)\cap H^m $ with respect to the norm $ \|.\|_{1,1} $, there must exist a $ \psi\in H^1 $ satisfying 
$ \langle\tilde{\varphi},\psi\rangle_{1,1}>0 $ and 
$$ \langle\varphi,\psi\rangle_{1,1}=0 \ \ \mbox{for all} \ \varphi\in{\rm Ker}(B)\cap H^m={\rm Ker}((A-1)^{-m}B)\cap H^m .$$
This last relation implies that for all $ \varphi\in{\rm Ker}((A-1)^{-m}B) $ one has
$$ \langle\varphi,(A-1)^{-m+1}\psi\rangle_{m,1}
       =\langle\varphi,(1-A)^{m-1}(A-1)^{-m+1}\psi\rangle_{1,1}
         =(-1)^{m-1}\langle\varphi,\psi\rangle_{1,1} = 0   , $$
which can be rephrased as
$$  (A-1)^{-m+1}\psi\perp_{m,1}{\rm Ker}((A-1)^{-m}B)   . $$
Since the operator $ (A-1)^{-m}B $ is anti-symmetric and compact with respect to the scalar-product 
$ \langle\cdot,\cdot\rangle_{m,1} $ it follows that
$$   (A-1)^{-m+1}\psi \in {\rm Range}((A-1)^{-m}B) .$$
Thus there exists a $ g\in H^m $ with $ (A-1)^{-1}Bg=\psi $ and it follows that
$$ \langle\tilde{\varphi},\psi\rangle_{1,1} = \langle\tilde{\varphi},(1-A)\psi\rangle
      =\langle\tilde{\varphi},(1-A)(A-1)^{-1}Bg\rangle = -\langle\tilde{\varphi},Bg\rangle=0 .$$
This contradicts our starting assumption. Therefore, all $ \tilde{\varphi}\in H^1 $ satisfying $ B\tilde{\varphi}\stackrel{w}{=}0 $ are in the closure of 
$ {\rm Ker}(B)\cap H^m $ with respect to the norm $ \|.\|_{1,1} $.
\hfill $ \Box $

\section{The limit is a resolvent again}

The operators $ A $ and $ A_c $ generate strongly continuous semigroups $ (T_t)_{t\geq0} $ resp. $ (T_t^{(c)})_{t\geq0} $
on $ H $. We will denote their respective resolvent families by $ (R_\lambda)_{\lambda>0} $ resp. $ (R_\lambda^{(c)})_{\lambda>0} $. 
From what we have seen so far follows that $ R_\lambda^{(c)}g=-\psi_c $ converges for $ |c|\rightarrow\infty $ toward 
$ R^*_\lambda g :=-\psi_* $ in $ H^1 $. We will see that once restricted to $ \overline{{\rm Ker}(B)\cap H^1}  $, 
the resulting limit is a strongly continuous resolvent, satisfying the contraction property.

\begin{theorem}\label{PseudoR}
The family $ (R_\lambda^*)_{\lambda>0} $ is a selfadjoint psudo-resolvent on $ H $; i.e. : 
$$ (\mu-\lambda)R_\lambda^\ast R_\mu^\ast=R_\lambda^\ast-R_\mu^\ast \ \ \ \mbox{for all $  \lambda,\mu>0$}. $$ 
\end{theorem}
{\bf Proof :}
The linearity is obvious.
The resolvent-identity holds for the true resolvent
$ (R_\lambda^{(c)})_{\lambda>0} $. One has for all $ g\in H $
$$  (\mu-\lambda)R_\lambda^{(c)}R_\mu^{(c)}g=R_\lambda^{(c)}g-R_\mu^{(c)}g  \ \ \ \mbox{for all $  \lambda,\mu>0$}.$$
Taking limits and using the fact, that for a fixed $ \lambda>0 $ the family $ \|R_\lambda^{(c)}\|;c\in\R $ is bounded,
one obtains the resolvent-identity for $ (R^*_\lambda)_{\lambda>0} $ :
$$  (\mu-\lambda)R_\lambda^*R_\mu^*=R_\lambda^*-R_\mu^* .$$
For $ \lambda>0 $ one has the representation (see Pazy (1983) p.25)
\begin{eqnarray*}
  R^{(c)}_\lambda f=\int_0^\infty e^{-\lambda t}T_t^{(c)}fdt.
\end{eqnarray*}
This implies 
\begin{eqnarray*}
 \langle R^{(c)}_\lambda f,g\rangle &=& \int_0^\infty e^{-\lambda t}\langle T_t^{(c)}f,g\rangle dt\\
 &=& \int_0^\infty e^{-\lambda t}\langle f, (T_t^{(c)})^*g\rangle dt\\
 &=& \int_0^\infty e^{-\lambda t}\langle f, T^{(-c)}_tg\rangle dt\\
 &=& \langle f,R_\lambda^{(-c)}g\rangle.
\end{eqnarray*}
Since by Theorem \ref{L5} both $ \psi_c=R^{(-c)}_\lambda g $ and $ \psi_{-c}=R^{(c)}_\lambda g $ converge to 
the same element $ R_\lambda^\ast g $ in $ H^1 $, one obtains in the limit
$$ \langle R^*_\lambda f,g\rangle=\langle f, R^*_\lambda g\rangle .$$
This shows that the family $ (R^\ast_\lambda)_{\lambda>0} $ is a self-adjoint pseudo-resolvent family (see Pazy (1983) p.36 for the definition of the 
pseudo-resolvent family).
\hfill $ \Box $

\vspace{5mm}

According to see Fukushima (1980) p.14 a resolvent family $ R_\lambda;\lambda>0 $ is said to have the contraction property, if 
$$  \langle\lambda R_\lambda g,\lambda R_\lambda g\rangle\leq\langle g,g\rangle \ \ \ \mbox{for all $ g\in H $ and $ \lambda>0 $.} $$
Following Fukushima (1980) p.23 we say that the resolvent family $ R_\lambda;\lambda>0 $ is Markovian if 
for all $ \lambda>0 $ and $ g\in H $ one has that
$$ 0\leq g\leq 1 \ \ \ \mbox{implies} \ \ \ 0\leq R_\lambda g\leq 1 .$$

\begin{corollary} \label{contraction}
The pseudo-resolvent family $ (R^\ast_\lambda)_{\lambda>0} $ has the contraction property and is Markovian.
\end{corollary}
{\bf Proof :}
It is well known that the resolvent family $ (R_\lambda^{(c)})_{\lambda>0} $ satisfies the contraction property.
Thus, we have for every $ \lambda>0 $, $ g\in H $ and $ c\in\R $ that 
$$   \left\langle\lambda R^{(c)}_\lambda g,\lambda R^{(c)}_\lambda g\right\rangle\leq \langle g,g\rangle .$$
Thus, since $ R^{(c)}_\lambda g $ converges in $ H $ toward $ R^\ast_\lambda g $ as $ c\rightarrow\infty $, we also have
$$   \left\langle\lambda R^\ast_\lambda g,\lambda R^\ast_\lambda g\right\rangle\leq \langle g,g\rangle .$$
Moreover, since the resolvent family $ (R_\lambda^{(c)})_{\lambda>0} $ is Markovian, we have for all $ \lambda>0 $, $ g\in H $ and $ c\in\R $ that 
$$ 0\leq g\leq 1 \ \mbox{$ \pi $-a.e.} \ \ \ \mbox{implies} \ \ \ 0\leq R^{(c)}_\lambda g\leq 1  \ \mbox{$ \pi $-a.e.}.$$
This implies for all $ c\in\R $ that
$$ \left\langle R^{(c)}_\lambda g,{\bf 1}_{R^\ast_\lambda g<0}\right\rangle \geq 0 .$$
From this follows 
$$ \left\langle R^\ast_\lambda g,{\bf 1}_{R^\ast_\lambda g<0}\right\rangle = 0 $$
and therefore $  R^\ast_\lambda g\geq 0 $ almost everywhere with respect to $ \pi $.\\[0.5mm]
On the other hand we also have 
$$ \left\langle 1-R^{(c)}_\lambda g,{\bf 1}_{1-R^\ast_\lambda g<0}\right\rangle \geq 0 .$$
From this follows 
$$ \left\langle 1-R^\ast_\lambda g,{\bf 1}_{1-R^\ast_\lambda g<0}\right\rangle = 0 $$
and therefore $  R^\ast_\lambda g\leq 1 $ almost everywhere with respect to $ \pi $.\\[0.5mm]
Thus we have proved the implication
$$ 0\leq g\leq 1 \  \mbox{$ \pi $-a.e.} \ \ \ \mbox{implies} \ \ \ 0\leq R^\ast_\lambda g\leq 1 \  \mbox{$ \pi $-a.e.} .$$
This finishes the proof of the proposition.
\hfill $ \Box $

\begin{proposition} \label{strongC}
The pseudo-resolvent family $ (R_\lambda^\ast)_{\lambda>0} $ is strongly continuous on the closure of 
$ {\rm Ker}(B)\cap H^1 $ in $ H $; i.e. : $$ \lim_{\lambda\rightarrow0}\|\lambda R_\lambda^*g-g\|\rightarrow0 \ \ \mbox{for all} \ 
 g\in\overline{{\rm Ker}(B)\cap H^1} .$$ 
\end{proposition}
{\bf Proof :}
First let $ g\in{\rm Ker}(B)\cap H^1 $. 
From the representation of  $ R^\ast_\lambda $ given in equation (\ref{ResolventenAusdruck}) we see that for $ g\in {\rm Ker}(B)\cap H^1 $ holds
\begin{eqnarray*}
\lambda R_\lambda^\ast g -g &=& -P_{1,\lambda}(\lambda(A-\lambda)^{-1}g+g) 
   = - P_{1,\lambda}(A(A-\lambda)^{-1}g).
\end{eqnarray*}
Now, recall that $ (\phi_n)_{n\in\N} $ is the ONB of $ H $ which is associated to the eigenvalues $ (-\lambda_n)_{n\in\N} $ of the selfadjoint operator $ A $; i.e. :
$$ A\phi_n=-\lambda_n\phi_n  . $$
We know that $ \lambda_n\geq0 $ for all $ n\in\N $ and $ \lambda_n\rightarrow\infty $ as $ n\uparrow\infty $. 
From the representation
$$ g=\sum_{n\in\N}\langle g,\phi_n\rangle\phi_n   $$
and the fact 
$$   \sum_{n\in\N}\lambda_n\langle g,\phi_n \rangle^2=\langle g,Ag\rangle<\infty   ,$$
we obtain for any $ g\in H^1 $ the following formula
$$ A(A-\lambda)^{-1}g=\sum_{n\in\N}\frac{-\lambda_n}{-\lambda_n-\lambda}\langle g,\phi_n\rangle\phi_n $$
This implies for $ g\in  H^1 $
\begin{eqnarray*}
                      \lambda\|A(A-\lambda)^{-1}g\|^2 
                 &=&   \sum_{n\in\N}\left(\frac{\lambda_n}{\lambda_n+\lambda}\right)^2\lambda\langle g,\phi_n\rangle^2 \\
                 &\leq& \sum_{n\in\N}\frac{\lambda_n^2}{\lambda_n+\lambda}\langle g,\phi_n\rangle^2   \longrightarrow 0 \ \ \ 
\mbox{ as $ \lambda\rightarrow\infty $. }
\end{eqnarray*}
Furthermore for $ g\in H^1 $ one has
$$  \|A(A-\lambda)^{-1}g\|_{H^1}^2=\sum_{n\in\N}\left(\frac{\lambda_n}{\lambda_n+\lambda}\right)^2\lambda_n\langle g,\phi_n\rangle^2 \longrightarrow0 \ \ \  \mbox{as $ \lambda\rightarrow\infty $.}  $$
Note that the convergence of the two above series follows from Lebesgue's theorem by using the majorant sequence
$ m_n :=\lambda_n\langle g,\phi_n\rangle^2 $ which satisfies
$$  \sum_{n\in\N}m_n=\sum_{n\in\N}\lambda_n\langle g,\phi_n\rangle^2=\|g\|_{H^1}^2<\infty .$$
Thus it follows for $ g\in {\rm Ker}(B)\cap H^1 $ that 
\begin{eqnarray*}
  \big\|\lambda R_\lambda^\ast g-g\big\|_{1,\lambda}^2 &=& \big\| P_{1,\lambda}(A(A-\lambda)^{-1}g)\big\|_{1,\lambda}^2 \\
       &\leq& \big\|A(A-\lambda)^{-1}g\big\|_{1,\lambda}^2 \\
     &=& \big\| A(A-\lambda)^{-1}g\big\|_{H^1}^2+\lambda\big\| A(A-\lambda)^{-1}g\big\|^2\longrightarrow0 \ \ \ \mbox{as $ \lambda\rightarrow\infty $}.
\end{eqnarray*}
In particular for $ g\in {\rm Ker}(B)\cap H^1 $ one has $ \|\lambda R_\lambda^\ast g-g\|\rightarrow0 $ as $ \lambda\rightarrow\infty $.
So far we have proved that the resolvent-family $ (R_\lambda^\ast) $ is strongly continuous on $ {\rm Ker}(B)\cap H^1 $.
Now let $ g $ be an element in the closure of $ {\rm Ker}(B)\cap H^1 $ and let $ \epsilon>0 $. 
Then there exists a $ \tilde{g}\in {\rm Ker}(B)\cap H^1 $ with the property $ \|\tilde{g}-g\|\leq\epsilon $.
Triangle inequality and the contraction property from Corollary \ref{contraction} then yield
\begin{eqnarray*}
 \|\lambda R^\ast_\lambda g - g\| &\leq& \|\lambda R^\ast_\lambda(g-\tilde{g})\|+\|\lambda R^\ast_\lambda\tilde{g}-\tilde{g}\|+\|\tilde{g}-g\| \\
 &\leq& \|\lambda R^\ast_\lambda\tilde{g}-\tilde{g}\|+2\|\tilde{g}-g\|.
\end{eqnarray*}
Thus the whole expression is smaller than $ 3\epsilon $ if $ \lambda $ is large enough. This proves the strong continuity on the closure of 
$ {\rm Ker}(B)\cap H^1 $ in $ H $.
\hfill $ \Box $

\begin{corollary}\label{ResolventeH1Konvergenz}
For any $ g\in {\rm Ker}(B)\cap H^1 $ one has that $ \lambda R_\lambda^\ast g $ converges in $ H^1 $ toward $ g $ as 
$ \lambda\rightarrow\infty $.
\end{corollary}
{\bf Proof :}
This follows from the proof of the previous proposition.
\hfill $ \Box $

\begin{proposition} \label{OrthoKern}
For all $ \lambda>0 $ we have that 
$$ {\rm Ker}(R^*_\lambda) =\left({\rm Ker}(B)\cap H^1\right)^\perp=\left(\overline{{\rm Ker}(B)\cap H^1}\right)^\perp .$$
Here, we use the definition
$$   ({\rm Ker}(B)\cap H^1)^\perp :=\left\{ g\in H; \langle g,h\rangle=0 \  \mbox{for all $ h\in{\rm Ker}(B)\cap H^1 $}\right\} .$$
\end{proposition}
{\bf Proof :} 
Since the second identity holds in general, we only have to prove the first.
We have for all $ g\in H $ and $ h\in{\rm Ker}(B)\cap H^1 $ that
$$   \langle g,h \rangle =-\langle (A-\lambda)^{-1}g,h\rangle_{1,\lambda}=-\langle P_{1,\lambda}(A-\lambda)^{-1}g,h\rangle_{1,\lambda}. $$
This shows that $ g\in\Big({\rm Ker}(B)\cap H^1\Big)^\perp $ is equivalent to $ P_{1,\lambda}(A-\lambda)^{-1}g=0 $.  
Since by (\ref{ResolventenAusdruck}) we have the representation $  R^\ast_\lambda g=-P_{1,\lambda}(A-\lambda)^{-1}g $, this proves the first 
identity. \hfill $ \Box $

\begin{proposition}
For all $ \lambda>0 $ we have that $ \overline{{\rm Range}(R^*_\lambda)}=\overline{{\rm Ker}(B)\cap H^1} $.
\end{proposition}
{\bf Proof :}
We will use the following representation for the range of a pseudo-resolvent, which can be found on page 217 in Yoshida 1980 :
$$ \overline{{\rm Range}(R_\lambda^*)}
=\left\{g\in H;\lim_{\lambda\rightarrow\infty}\lambda R_\lambda^*g=g\ {\rm in}\ H \right\} .$$
Proposition \ref{strongC} shows that for all $ g\in\overline{{\rm Ker}(B)\cap H^1} $ we have $ \lambda R_\lambda^*g\rightarrow g $ as $ \lambda\rightarrow\infty $. The inclusion $  \overline{{\rm Ker}(B)\cap H^1}\subset\overline{{\rm Range}(R_\lambda^*)} $ then follows.\\[0.5mm]
To prove the inclusion $  \overline{{\rm Range}(R_\lambda^*)}\subset\overline{{\rm Ker}(B)\cap H^1} $ it is enough to show the inclusion $  {\rm Range}(R_\lambda^*)\subset{\rm Ker}(B)\cap H^1 $. This however follows from the fact that we have the representation
$ R_\lambda^\ast g=-P_{1,\lambda} (A-\lambda)^{-1}g $.
\hfill $ \Box $

\begin{proposition}
For all $ g\in{\rm Range}(B) $ and $ \lambda>0 $ one has $ R_\lambda^{(c)}g\rightarrow0 $ as $ |c|\rightarrow\infty $.
\end{proposition}
{\bf Proof:} As $ B $ is an antisymetric operator, one has that the subspaces $ {\rm Ker}(B) $ and $ {\rm Range}(B) $ are orthogonal. Any $ h\in H $ can be decomposed as $ h=h_1+h_2 $ where $ h_1\in{\rm Ker}(B) $ and $ h_2\in{\rm Range}(B) $. It then follows that
\begin{eqnarray*}
   \langle h,R_\lambda^\ast g\rangle &=& \langle h_1,R_\lambda^\ast g\rangle \ \ \ \mbox{since $ R_\lambda^\ast g\in{\rm Ker}(B) $} \\
    &=& \langle R_\lambda^\ast h_1,g\rangle \ \ \ \mbox{as $ R_\lambda^\ast $ is self-adjoint} \\
    &=&  0 \ \ \ \mbox{since $ {\rm Ker}(B)\perp{\rm Range}(B) $}.
\end{eqnarray*}
Since this holds for all $ H\in H $, it follows that  $ R_\lambda^\ast g=0 $.
\hfill $ \Box $

\section{The quest for a suitable state-space} \label{Sektion5}

In order to associate the resolvent family $ (R_\lambda^\ast)_{\lambda>0} $ to a diffusion, we have to prove, that its domain of definition 
 $ \overline{{\rm Ker}(B)\cap H^1} $ can be identified with some space of square integrable functions over some suitable 
measure space $ (X,\mu) $. In this section we give a construction for this space. 
The vectorfield $ b :M\rightarrow TM $ generates a flow $ (\Phi_t)_{t\in\mathbb{R}} $ on $ M $.
We say that two elements $ x $, $ y $ from $ M $ are equivalent, if there exist for all $ \epsilon>0 $ a $ N\in\N $, a sequence $ (z_n)_{n=0,...,N} $ in $ M $
with $ z_0=x $, $ z_N=y $ and two sequences of real numbers $  (t_n)_{n=0,...,N-1} $ and $ (s_n)_{n=1,...,N} $ with the property
$$      d(\Phi_{t_0}(z_0),\Phi_{s_1}(z_1))+ d(\Phi_{t_1}(z_1),\Phi_{s_2}(z_2))+...+ d(\Phi_{t_{N-1}}(z_{N-1}),\Phi_{s_N}(z_N))  < \epsilon .$$
We will note $ x\sim y $ when two elements $ x $, $ y $ from $ M $ are equivalent. 
For an $ x\in M $ we denote by $ [x] $ the equivalence class containing $ x $. 
Let $ X $ be the quotient space with respect to this equivalence relation. 
The space $ X $ carries the quotient topology $ \mathcal{O}_q $; i.e. : the largest topology on $ X $ with the property that the projection $ p_X :M\rightarrow X;x\mapsto[x] $ is continuous.

\begin{proposition}
The quotient topology on $ X $ can be obtained from the following metric
$$ d_X\big([x],[y]\big)  :=\inf\left\{\sum_{i=0}^{N-1} d(\Phi_{t_i}(z_i),\Phi_{s_{i+1}}(z_{i+1})); N, (t_n),(s_n), (z_n) :  z_0=x,z_N=y \right\}. $$
\end{proposition}
{\bf Proof:}
One has to prove that the topology $ \mathcal{O}_d $ generated by the metric $ d_X $ equals $ \mathcal{O}_q $. 
In order to show the inclusion $ \mathcal{O}_d\subset\mathcal{O}_q $, we have to prove that the projection $ p_X:M\rightarrow X;x\mapsto[x] $ is continuous with respect to the metric $ d_X $. This follows from the following inequalities
\begin{eqnarray*}
 d_X([x],[y]) \!\!\! &=& \!\!\! \inf\left\{\sum_{i=0}^{N-1} d(\Phi_{t_i}(z_i),\Phi_{s_{i+1}}(z_{i+1})); N, (t_n),(s_n), (z_n) :  z_0=x,z_N=y \right\} \\
 \!\!\!   &\leq& \!\!\! d(x,y) .
\end{eqnarray*}
In order to prove the inclusion $ \mathcal{O}_q\subset\mathcal{O}_d $, we assume that there exist some $ U\in\mathcal{O}_q $ which is not in $ \mathcal{O}_d $.
This means that $ U $ is not open with respect to $ d_X $. 
It follows that there exists some $ [x]\in U $ and a sequence $ [y_n]\in U^c $ such that $ d_X([x],[y_n])\rightarrow0 $ as $ n\rightarrow\infty $.
We then have that $ \tilde{U}=p_X^{-1}(U) $ is open in $ M $ and $ x'\in\tilde{U} $ for all $ x'\in[x] $. 
Further, there then exists for all $ x'\in[x] $ some $ \epsilon>0 $ such that $ d(x',y_n')>\epsilon $ for all selections of $ y'_n\in[y_n] $ and $ n\in\N $. 
If we consider one of those subsequences $ y'_n;n\in\N $, there then exists a subsequence $ y'_{n_k};k\in\N $ which converges toward a limit  point $ x'' $ in the compact manifold $ M $. 
It then follows that for all $ \epsilon'>0 $ there exists some $ K\in\N $ such that $ d(y_{n_k}',x'')<\epsilon' $ and $ d_X([y_{n_k}'],[x])<\epsilon' $ for all $ k\geq K $. 
From this follows that $ x''\in[x] $ and $ y_{n_k}'\rightarrow x'' $ as $ k\rightarrow\infty $. 
This is a contradiction and thus there is no $ U\in\mathcal{O}_q $ which is not also in the topology $ \mathcal{O}_d $.
\hfill   $ \Box $

\bigskip

With this metric the space $ (X,d) $ is a compact Hausdorf space. 
We denote by $ {\mathcal B}(X) $ the Borel-$ \sigma $-field generated by this topology on $ X $. 
Then the map $ p_X $ is measurable with respect to the respective Borel-$ \sigma $-fields on $ M $ and $ X $.
On $ M $ we define the eventually smaller $ \sigma $-field 
$$ \Sigma_M :=\left\{p^{-1}_X(A);A\in {\mathcal B}(X)\right\}  .$$
We introduce the measure $ \mu :=\pi\circ p_X^{-1} $ on $ (X,{\mathcal B}(X)) $ and note that since $ \pi $ has full support on $ M $ the measure $ \mu $ also has full support on $ X $.\\[0.5mm]
For a $ \varphi\in L^2(M,\pi) $ we can define the orthogonal projection from $ L^2(M,\pi) $ to the closed subspace 
$ I^2(M) :=\{f \in L^2(M,\pi) :\ \mbox{$ f $ is $ \Sigma_M $-measurable}\} $.
We denote this projection by $ P_I $ and note that by factorization lemma for every $ \varphi\in L^2(M,\pi) $ which is 
$ \Sigma_M $-measurable there exists a unique $ \tilde{\varphi}\in L^2(X,\mu) $ with the property 
$ \varphi=\tilde{\varphi}\circ p_X $. The above consideration defines a map
$$ R:L^2(M,\Sigma_M,\pi)\rightarrow L^2(X,\mu);\varphi\mapsto\tilde{\varphi} .$$
We define the following two maps :
\begin{eqnarray} \label{GIsometry}
 G :\overline{{\rm Ker}(B)\cap H^1}\rightarrow L^2(X,\mu); \varphi\mapsto R(P_I(\varphi)) 
\end{eqnarray}
and 
\begin{eqnarray} \label{FIsometry}
 F :L^2(X,\mu)\rightarrow\overline{{\rm Ker}(B)\cap H^1}; \psi\mapsto\Big[x\mapsto \psi([x])\Big] .
 \end{eqnarray}
Note that although the map $ G $ could be defined on all of $ L^2(M,\pi) $, we restrict its domain to $ \overline{{\rm Ker}(B)\cap H^1} $.
We have to prove that $ G $ and $ F $ are well defined. Obviously $ G $ is well defined, since the involved maps $ P_I $ and $ R $ 
are well defined. In order to prove that $ F $ is well defined we have to see that the resulting function $ x\mapsto \psi([x]) $ 
is indeed an element from $ \overline{{\rm Ker}(B)\cap H^1} $. We know from Theorem \ref{regularity} that this space is equal to 
$ \{\varphi\in H :B\varphi\stackrel{w}{=}0\} $. Thus we have to see whether $ \tilde{\psi}(x) :=\psi([x]) $ satisfies 
$ \langle \tilde{\psi},B\vartheta\rangle=0 $ for all $ \vartheta\in H^1 $. We have for all $ t\in\mathbb{R} $ that $ [\Phi_t(x)]=[x] $ and therefore
$$    \langle \tilde{\psi},\vartheta\circ\Phi_t\rangle = \langle \tilde{\psi}\circ\Phi_{-t},\vartheta\rangle
                                                       = \langle \tilde{\psi},\vartheta\rangle                   .$$
Thus we have
$$  \langle \tilde{\psi},B\vartheta\rangle=
 \lim_{t\rightarrow0}\frac{1}{t}\Big(\langle \tilde{\psi},\vartheta\circ\Phi_t\rangle- \langle \tilde{\psi},\vartheta\rangle  \Big)=0 .$$
It therefore turns out that $ F $ is well defined and also continuous, since it is linear and bounded. \\[1ex]

\begin{lemma} \label{levelsets}
For a Lipshitz continuous $ \varphi\in C(M)\cap{\rm Ker}(B) $ one has that $ x\sim y $ implies $ \varphi(x)=\varphi(y) $. As a consequence, the assignment  $ \tilde{\phi}([x]):=\phi(x) $ defines a well defined and continuous function on $ X $.
\end{lemma}
{\bf Proof :}
Let $ \varphi\in  C(M)\cap{\rm Ker}(B) $. 
According to Corollary \ref{Dicht}, we can find some sequence $ \varphi_n\in C^1(M)\cap {\rm Ker}(B) $ such that $ \varphi_n\rightarrow\varphi $ in $ H $.
First note that for all $ n\in\mathbb{N} $ one has $ 0=B\varphi_n=b\cdot\nabla\varphi_n=\partial_t\big|_{t=0}\varphi_n\circ\Phi_t $. 
This implies that one has $ \varphi_n\circ\Phi_t(x)=\varphi_n(x) $ for all $ x\in M, t\in\mathbb{R} $ and all  $ n\in\mathbb{N} $. 
It then follows that some sub-sequence $ \varphi_{n(k)} $ converges toward $ \varphi $ almost everywhere. 
Taking limit $ n\rightarrow\infty $ this implies that for any choice of $ t\in\mathbb{R} $ one has $ \varphi\circ\Phi_t(x) = \varphi(x) $ almost everywhere with respect to $ x $. 
Since $ \varphi $ and $ \varphi\circ\Phi_t $ are continuous, it follows that  $ \varphi\circ\Phi_t(x) = \varphi(x) $ for all $ x\in M, t\in\mathbb{R} $. 
The  compactness of $ M $ implies that $ \varphi $ is uniformly continuous on $ M $; i.e. : for all $ \epsilon>0 $ there exist a $ \delta>0 $ such that 
$$ d(x,y)<\delta \ \ \ \mbox{implies} \ \ \ \Big|\varphi(x)-\varphi(y)\Big|<\epsilon .$$
Let $ \eta>0 $. Since $ x\sim y $ there must exist $ N $, $ (t_n) $, $ (s_n) $ and $ (z_n) $ with $ z_0=x $ and $ z_N=y $ such that
$$ \sum_{i=0}^Nd(\Phi_{t_i}(z_i),\Phi_{s_i}(z_{i+1}))< \eta .$$
It follows from this and the fact that $ \varphi $ is constant along the flow that
\begin{eqnarray*}
   \Big|\varphi(x)-\varphi(y)\Big| &\leq& \sum_{i=0}^N\Big|\varphi\big(\Phi_{t_i}(z_i)\big)-\varphi\big(\Phi_{s_i}(z_{i+1})\big)\Big| \\
    &\leq& C\sum_{i=0}^N\Big|d\big(\Phi_{t_i}(z_i),\Phi_{s_i}(z_{i+1})\big)\Big| \ \leq \ \epsilon ,
\end{eqnarray*}
where $ C>0 $ is the Lipschitz constant of the function $ \phi $.
Since $ \epsilon>0 $ can be arbitrarily small this proves the statement. 
\hfill $ \Box $

\vspace{1ex}

{\bf Remark :} The previous lemma can't be generalized to general continuous functions. A counter example showing that this is not possible can be found in Section \ref{Beispiele} (see example 6). 

\begin{theorem} \label{Isometrie}
The maps $ F $ and $ G $ are bijective and $ G $ is the inverse map of $ F $. Moreover, $ G $ is an isometry between the spaces 
$ \overline{{\rm Ker}(B)\cap H^1} $ and $ L^2(X,\mu) $.
\end{theorem}
{\bf Proof :}
We first prove that $ G\circ F $ is the identity map on $ L^2(X,\mu) $. This then implies that $ F $ is injective and $ G $ 
is surjective. For a $ \psi\in L^2(X,\mu) $ the function $ F(\psi)=\psi\circ p_X $ is measurable with respect to $ \Sigma_M $.  
Thus we have $ P_I(F(\psi))=\psi\circ p_X $. This then implies that $ G\circ F(\psi)=R(\psi\circ p_X)=\psi $.\\[0.5mm]
It remains to prove that the map $ F $ is surjective and the map $ G $ is injective. \\[0.5mm]
Since the map $ G $ is linear, in order to prove its injectivity, it is sufficient to prove that its kernel only 
contains the zero-function from $ \overline{{\rm Ker}(B)\cap H^1} $. 
Let $ \varphi\in\overline{{\rm Ker}(B)\cap H^1} $ be some element with $ G(\varphi)=0 $. This means that $ P_I(\varphi)=0 $. 
We know from Corollary \ref{Dicht} that there exists a sequence of $ C^1 $-functions 
$ \varphi_n\in{\rm Ker}(B)\cap H^1 $ which converges toward $ \varphi $ in $ H $. 
Since those elements are continuously differentiable, they satisfy $$ 0=B\varphi_n=b\cdot\nabla\varphi_n .$$ 
This implies that $ \varphi_n $ is invariant with respect to the flow $ (\Phi_t)_{t\in\R} $ in the strong sense; i.e. : $ \varphi_n\circ\Phi_t(x)=\varphi_n(x) $ for all $ t\in\R $ and all $ x\in M $. 
From Lemma \ref{levelsets} we know that the values of $ \varphi_n $ are constant on the equivalence classes and it then follows that one has $ P_I(\varphi_n)=\varphi_n $ for all $ n\in\N $.
Moreover, we know that in $ L^2(M,\pi) $ we have $ P_I(\varphi_n)\rightarrow P_I(\varphi)=0 $. Therefore, $ \varphi_n\rightarrow0 $.
Since $ \varphi_n\rightarrow\varphi $, we conclude $ \varphi=0 $, which proves the injectivity of $ G $.\\[0.5mm]
It remains to prove the surjectivity of $ F $. For a $ \varphi\in \overline{{\rm Ker}(B)\cap H^1} $ we saw in the first part of this proof 
that there exists a sequence of $ C^1 $-functions $ \varphi_n $ which are constant on the equivalence classes defined on $ M $ with
$ \varphi_n\rightarrow\varphi $ in $ L^2(M,\pi) $. Obviously there exist functions $ \psi_n\in L^2(X,\mu) $ with the property that 
$ \psi_n([x])=\varphi_n(x) $ for all $ x\in M $; i.e. : $ F(\psi_n)=\varphi_n $. Moreover, by the definition of $ \mu $ we see that
\begin{eqnarray*}
 \int_X|\psi_n-\psi_m|^2d\mu 
    &=&    \int_M|\psi_n\circ p_X-\psi_m\circ p_X|^2d\pi \\
       &=&  \int_M|\varphi_n-\varphi_m|^2d\pi .
 \end{eqnarray*}
The convergence of $ \varphi_n $ in $ L^2(M,\pi) $ then shows, that the sequence $ \psi_n; n\in\N $ is a Cauchy-sequence and should converge 
toward some element $ \psi $ in $ L^2(X,\mu) $. The continuity of $ F $ then implies that $ F(\psi)=\varphi $.\\[0.5mm]
Since, $ F $ and $ G $ are both bijective and $ G\circ F $ is the identity on $ \overline{{\rm Ker}(B)\cap H^1} $ it follows that $ G $ is the inverse map of $ F $.
\hfill $ \Box $

\vspace{5mm}

Note that by the previous definition for an element $ h\in L^2(M,\pi) $ the image $ G(h) $ is an element in $ L^2(X,\mu) $, which is formally an equivalence class of indistinguishable functions. 
This means that there is some ambiguity on the choice of some representing element from the class. However, in the case of a continuous function $ h $, which are constant on the equivalence classes, there might be a canonical choice such that $ G(h) $ is also continuous. In order to state this result, we introduce the space
$$ C^\sim(M)=\Big\{f\in C(M); f(x)=f(y) \ \mbox{whenever} \ x\sim y  \Big\} .$$

\begin{proposition} \label{CBild}
For each $ h\in C^\sim(M)\cap\overline{{\rm Ker}(B)\cap H^1} $, there is a unique element $ \tilde{h} $ in the equivalence class $ G(h) $ satisfying $ \tilde{h}\in C(X) $.
If for $ h\in C(M) $ we identify the class $ G(h) $ with $ \tilde{h} $, then one has 
$$ C(X)=G\Big(C^\sim(M)\cap\overline{{\rm Ker}(B)\cap H^1}\Big) .$$
\end{proposition}
{\bf Proof :} Let $ h\in C(X) $ and $ (y_n)_{n\in\N} $ a sequence in $ M $ which converges toward some $ y\in M $. Since we always have 
$ d_X([u],[v])\leq d_M(u,v) $, the sequence $ ([y_n])_{n\in\N} $ converges toward $ [y] $ in $ X $. 
Thus $ F(h)(y_n)=h([y_n]) $ converges toward $ h([y])=F(h)(y)$, which proves the continuity of $ F(h) $. By construction $ F(h) $ is then also in $ C^\sim(M) $.
This implies $$ F(C(X))\subset C^\sim(M)\cap\overline{{\rm Ker}(B)\cap H^1}$$ from which follows by Theorem \ref{Isometrie} that
$$ C(X)\subset G\Big(C^\sim(M)\cap\overline{{\rm Ker}(B)\cap H^1}\Big)  .$$
We still have to prove that $  G\Big(C^\sim(M)\cap\overline{{\rm Ker}(B)\cap H^1}\Big) \subset C(X) $. For this let $ h\in C^\sim(M)\cap\overline{{\rm Ker}(B)\cap H^1} $.From this follows that $ h(x)=h(y) $ when $ [x]=[y] $. 
Therefore, the function $ \tilde{h}([x])=h(x) $ is well defined on $ X $ and one has $ h=\tilde{h}\circ p_X $.
It follows from the properties of the quotient topology on $ X $, that if $ h $ is continuous on $ M $ then $ \tilde{h} $ has to be continuous on $ X $.
As a consequence the function $ \tilde{h} $ is measrable with respect to the Borel-$ \sigma $-field $ \mathcal{B}(X) $. 
This implies that $ h $ is $ \Sigma_M $-measurable and thus an element from $ I^2(M) $. 
From this follows that $ P_I(h)=h $, which yields
$$  G(h)=R(P_I(h))=R(h)=\tilde{h}, \ \mbox{a.e.} . $$ 
The abovementioned continuity of $ \tilde{h} $ then shows that $ G(h)\in C(X) $.
\hfill $ \Box $

{\begin{corollary}
One has $ I^2(M)=\overline{{\rm Ker}(B)\cap H^1} $. 
\end{corollary}
{\bf Proof:} We have seen in the previous proof that $ \mathcal{C}^\sim(M)\cap\overline{{\rm Ker}(B)\cap H^1}\subset I^2(M) $
and we know that $ C^\sim(M)\cap\overline{{\rm Ker}(B)\cap H^1} $ is dense in $ \overline{{\rm Ker}(B)\cap H^1} $.
This shows that $ I^2(M)\subset\overline{{\rm Ker}(B)\cap H^1} $. 
On the other hand, each element $ \psi\in\overline{{\rm Ker}(B)\cap H^1}$ is a limit in $ H $ of $ C^1 $-functions $ \psi_n $ from $ {\rm Ker}(B) $. 
Since $ \psi_n\in I^2(M) $ for all $ n\in\N $, it follows that their limit $ \psi $ is also in $ I^2(M) $. 
This shows that $ \overline{{\rm Ker}(B)\cap H^1}\subset I^2(M) $ and thus $ I^2(M)=\overline{{\rm Ker}(B)\cap H^1} $.
\hfill $ \Box $

\begin{proposition} \label{Homeomorphism}
The map $ G $ is an algebra homeomorphism between the two algebras 
$ C^\sim(M)\cap\overline{{\rm Ker}(B)\cap H^1} $ and $ C(X) $.
\end{proposition}
{\bf Proof :}
It is obvious from the definition, that $ F $ is an algebra homeomorphism between $ C(X) $ and $ F(C(X)) $ which equals $ C^\sim(M)\cap\overline{{\rm Ker}(B)\cap H^1} $ by Proposition \ref{CBild}. From Theorem \ref{Isometrie} we know that the map $ G  $ is the inverse of $ F $. It follows that $ G $ has to be an algebra homeomorphism from $ C^\sim(M)\cap\overline{{\rm Ker}(B)\cap H^1} $ to $ C(X) $.
\hfill $ \Box $

\section{The generator and the Dirichlet form}

In this section we compute the generator and the Dirichlet-form underlying the limit resolvent $ (R_\lambda^\ast)_{\lambda>0} $. 

\begin{theorem} \label{Dform}
The Dirichlet-form $ {\mathcal E}^* $ associated to the resolvent family $ (R_\lambda^\ast)_{\lambda>0} $  is the closed form defined on 
$ {\rm Ker}(B)\cap H^1 $ as follows
$$   {\mathcal E}^* : (u,v)\mapsto\langle u, v\rangle_{H^1}  .$$
\end{theorem}
{\bf Proof :}
We have for $ u,v\in{\rm Ker}(B)\cap H^1 $ :
\begin{eqnarray*}
  \lambda\big\langle u-\lambda R^\ast_\lambda u,v\big\rangle 
&=& \lambda\big\langle u+\lambda P_{1,\lambda}(A-\lambda)^ {-1}u,v\big\rangle \\
&=& \lambda\Big(\big\langle u,v\big\rangle +\big\langle P_{1,\lambda}(A-\lambda)^{-1}u,v\big\rangle_{1,\lambda}-\big\langle P_{1,\lambda}(A-\lambda)^{-1}u,v\big\rangle_{H^1}\Big)\\
&=& \lambda\Big(\big\langle u,v\big\rangle -\big\langle u,P_{1,\lambda}v\big\rangle -\big\langle P_{1,\lambda}(A-\lambda)^{-1}u,v\big\rangle_{H^1}\Big) \\
&=&  -\lambda\big\langle P_{1,\lambda}(\lambda-A)^{-1} u,v\big\rangle_{H^1}   \ \ \ \ \ \ \ \ \ \ \ \ \ \ \ \ \ \ \ \ \ \ \ \ \mbox{(since $ P_{1,\lambda}v=v $)}\\
&=& \big\langle\lambda R^\ast_\lambda u,v\big\rangle_{H^1} .
\end{eqnarray*}
Corollary \ref{ResolventeH1Konvergenz} shows that the last term converges for $ \lambda\rightarrow\infty $ toward
$$ \big\langle u,v\big\rangle_{H^1}={\mathcal E}^*(u,v) .$$ 
This characterizes the associated Dirichlet form (see Fukushima (1980) p.21).\\[0.5mm]
To prove that the form is closed, we have to prove that for every sequence 
$ \psi_n\in{\rm Ker}(B)\cap H^1 $ satisfying 
$  {\mathcal E}^*(\psi_n-\psi_m,\psi_n-\psi_m)\rightarrow0 $ as $ n,m\uparrow\infty $ there exists some $ \psi\in{\rm Ker}(B)\cap H^1 $ with
$ {\mathcal E}^*(\psi_n-\psi,\psi_n-\psi)\rightarrow0  $ (see Fukushima p.4).
Since we know that $  (\varphi,\psi)\mapsto\langle\varphi,\psi\rangle_{H^1} $ is a closed form on $ H^1 $, we know that there exists some 
$ \psi\in H^1 $ with $ \langle\psi_n-\psi,\psi_n-\psi\rangle_{H^1}\rightarrow0  $.
We have to prove that $ \psi $ is also in $ {\rm Ker}(B) $.
This follows from the fact that $ {\rm Ker}(B) $ is closed in $ H^1 $.
\hfill $ \Box $\\

We now obtain a Dirichlet form on the space $ L^2(X,\mu) $ through the maps $ G $ and $ F$ from Theorem \ref{Isometrie}. It will turn out that the resulting Dirichlet form is regular (see Fukushima p.6 for a definition).
We saw in Proposition \ref{CBild} that for any function $ h $  in the set $ S=C^\sim(M)\cap{\rm Ker}(B)\cap H^1 $ there exist some unique function $ \tilde{h} $ in the equivalence class $ G(h) $, which is continuous on $ X $.
For the statement of the following proposition, we identify for $ h\in S $ the class $ G(h) $ with the function $ \tilde{h} $.

\begin{proposition} \label{Cdicht}
The set $ G(S) $ is dense in $ C(X) $ with respect to sup-norm.
\end{proposition}
{\bf Proof :} We first note that $ C^\sim(M)\cap{\rm Ker}(B)\cap H^1 $ is a commutative algebra. In order to use the Weierstrass theorem(see Dieudonn\'e 1969) it remains to show that $ G(S) $ contains the constant functions and that for two different points $ [x],[y]\in X $ there exists a $ \xi\in S $ with $ G(\xi)([x])\neq G(\xi)([y]) $. 
The first part is obvious and so we have to deal with the second one. Let $ d_M(.,.) $ and $ d_X(.,.) $ denote the metrics on $ M $ and $ X $ respectively. 
Obviously the function $ h :x\mapsto d_ X([x],[y]) $ is in $ {\rm Ker}(B) $ and seperates the two points $ [x] $ and $ [y] $. It is sufficient to show that $ h $ is an element from $ H^1 $. 
For this we note that 
$$ d_X([x],[y])-d_X([x'],[y])\leq d_X([x],[x']) \leq d_M(x,x') \ \ \  \forall  x,x'\in M .$$ 
This proves that the function $ h $ is Lipschitz continuous. Now, by Rademacher theorem this implies that  $ h $ is in $ H^1 $. 
\hfill $ \Box $

\vspace{5mm}

Let $ m $ be a measure on the Borel $ \sigma $-field of a compact metric space $ (Y,d) $. 
According to Fukushima (1980, p.5 and 6) a set $ C_1\subset\mathcal{D}(\mathcal{E})\cap C(Y) $ is said to be a core for a Dirichlet form $ (\mathcal{E},\mathcal{D}(\mathcal{E})) $ on $ L^2(Y,m) $ if the two following properties hold :\\[0.5mm]
i) $ C_1 $ is dense in $ \mathcal{D}(\mathcal{E}) $ with respect to the norm $ \|u\|^2_1:=\mathcal{E}(u,u)+\|u\|^2_{L^2(Y,m)} $;\\[0.5mm]
ii) $ C_1 $ is dense in $ C(Y) $ with respect to the supremum norm $ \|.\|_{\sup} $.\\[0.5mm]
A Dirichlet form $ (\mathcal{E},\mathcal{D}(\mathcal{E})) $ on $ L^2(Y,m) $ is called regular if it has a core.

\begin{corollary} 
The form $ {\mathcal E}_*(\varphi,\psi) :={\mathcal E}^*(F(\varphi),F(\psi)) $  with domain $  G( {\rm Ker}(B)\cap H^1 ) $ is a regular Dirichlet form on $ L^2(X,\mu) $. 
\end{corollary}
{\bf Proof :}
We have to show that 
$ C(X)\cap  G( {\rm Ker}(B)\cap H^1) $ is dense in $  G( {\rm Ker}(B)\cap H^1 ) $ with respect to the norm defined on $  G( {\rm Ker}(B)\cap H^1 ) $ through
$$ \|\psi\|_{1,1,X}^2 :={\mathcal E}_*(\psi,\psi)+\|\psi\|_{L^2(X,\mu)}^2 $$ 
and dense in $ C(X) $ with respect to the supremum norm. 
By Sobolev embedding theorem and Proposition \ref{H1dicht} it follows that $ C^1(M)\cap{\rm Ker}(B)\cap H^1 $ is dense in $ {\rm Ker}(B)\cap H^1 $ with respect to the $ \|.\|_{1,1} $-norm. Since the functions from $ C^1(M)\cap{\rm Ker}(B)\cap H^1 $ are constant on equivalence classes by Lemma \ref{levelsets}, it follows that  $ C^\sim(M)\cap{\rm Ker}(B)\cap H^1 $ is dense in $ {\rm Ker}(B)\cap H^1 $ with respect to the same norm.
If $ \psi\in G( {\rm Ker}(B)\cap H^1 ) $ there exists some $ \varphi\in{\rm Ker}(B)\cap H^1 $ such that $ G(\varphi)=\psi $ and thus by Theorem 4 one has $ F(\psi)=\varphi $. 
As we just saw, there exists a sequence $ \varphi_n;n\in\mathbb{N} $ from $  C^\sim(M)\cap{\rm Ker}(B)\cap H^1 $ satisfying 
$$   {\mathcal E}^*(\varphi_n-\varphi,\varphi_n-\varphi)+\|\varphi_n-\varphi\|^2=\|\varphi_n-\varphi\|_{1,1}^2\longrightarrow 0 \ \mbox{as $ n\rightarrow\infty $}. $$
From the fact that $ G $ is an isometry from $ \overline{{\rm Ker}(B)\cap H^1} $ to $ L^2(X,\mu) $ and by the defintion of the bilinear form $ {\mathcal E}_* $ we obtain from this
\begin{eqnarray*} 
  {\mathcal E}_*(G(\varphi_n)-\psi,G(\varphi_n)-\psi)+\|G(\varphi_n)-\psi\|^2_{L^2(X,\mu)}\longrightarrow 0 \ \mbox{as $ n\rightarrow\infty $}. 
\end{eqnarray*}
Since we know from Proposition \ref{CBild} that the elements of the sequence $ G(\varphi_n);n\in\mathbb{N} $ can be assumed to be contiuous, we have proved that 
$ C(X)\cap  G( {\rm Ker}(B)\cap H^1 ) $ is dense in $  G( {\rm Ker}(B)\cap H^1 ) $  with respect to the norm 
$ \|.\|_{1,1,X} $ introduced previously.
In order to prove the second requirement, let $ \tilde{\psi}\in C(X) $. 
By Proposition \ref{Cdicht} there exists a sequence  $ \tilde{\varphi}_n;n\in\mathbb{N} $ from
 $  C^\sim(M)\cap{\rm Ker}(B)\cap H^1 $ satisfying
$$  \|G(\tilde{\varphi}_n)-\tilde{\psi}\|_{\infty}\longrightarrow0 \  \mbox{as $ n\rightarrow\infty $}, $$
where $ \|.\|_\infty $ denotes the sup-norm on $ C(X) $. 
Clearly $ \tilde{\psi}_n :=G(\tilde{\varphi}_n) $ is in $ G({\rm Ker}(B)\cap H^1 ) $ and further one has by Proposition \ref{CBild} that $\tilde{\psi}_n\in C(X) $.
This shows that the set $ C(X)\cap  G( {\rm Ker}(B)\cap H^1 ) $ is dense in $ C(X) $ with respect to $ \sup $-norm.
\hfill $ \Box $

\vspace{5mm}

The Dirichlet form $ {\mathcal E}^* $ is associated to some strongly continuous contraction semigroup $ (T^\ast_t)_{t\geq0} $ of symmetric operators 
on $ \overline{{\rm Ker}(B)\cap H^1}  $ (see Fukushima p.14ff). Moreover, the semigroup  $ (T^\ast_t)_{t\geq0} $ is related to the resolvent family 
$ (R_\lambda^\ast)_{\lambda>0} $ through 
$$   R_\lambda^\ast g=\int_0^\infty e^{-\lambda t}T^\ast_t gdt .$$
According to Fukushima (1980) p.15 there exists a generator 
$$ A^\ast :{\mathcal D}(A^\ast)\rightarrow H; f\mapsto\lim_{t\rightarrow0}\frac{T_t^\ast f-f}{t} $$ 
with domain of definition
$$ {\mathcal D}(A^\ast) :=\left\{f\in H;\lim_{t\rightarrow0}\frac{T_t^\ast f-f}{t} \ \mbox{exists in $ H $}\right\}   .$$
According to general semi-group theory (see Pazy (183)) one has that :\\[0.5mm]
i) $ \int_0^tT^\ast_s\varphi ds\in {\mathcal D}(A^\ast) $ for all $ t>0 $ and $ \varphi\in\overline{{\rm Ker}(B)\cap H^1} $;\\[0.5mm]
ii) $ \frac{d}{dt}T_t^\ast\varphi=T_t^\ast A^\ast\varphi  $ for all $ \varphi\in{\mathcal D}(A^\ast) $;\\[0.5mm]
iii) $ {\rm Range}(R_\lambda^\ast)={\mathcal D}(A^\ast) $.\\[0.5mm]

In order to obtain some more explicit representation for the generator $ A^\ast $ we will use $ P_E $ to denote the orthogonal projection from $ H $ to the sub-space  $ E :=\overline{{\rm Ker}(B)\cap H^1} $. We have the following result.

\begin{proposition} \label{Erzeuger}
Restricted to the space $ E $ the family $ (R^*_\lambda)_{\lambda>0} $ 
is the resolvent of the closure of the operator $ A^* $ which on  $ {\rm Ker}(B)\cap H^2  $ has the representation
$$ A^\ast g=P_E Ag .$$
\end{proposition}
{\bf Proof :}
Note that we have the representation $ R_\lambda^\ast =-P_{1,\lambda}(A-\lambda)^{-1}  $.
In the following argument, we denote $ C' $ the adjoint of a bounded operator $ C $. 
Since $ R_\lambda^\ast $ is self-adjoint and since its range is a subset of $ E $, we have that 
$$ R_\lambda^\ast =(R_\lambda^\ast)' =(P_E R_\lambda^\ast)'=(R_\lambda^\ast)' P_E =R_\lambda^\ast P_E .$$
It then follows for any $ g\in {\rm Ker}(B)\cap H^2 $ that 
\begin{eqnarray*}
  R_\lambda^\ast (\lambda-P_EA)g &=&
  R_\lambda^\ast P_E(\lambda-A)g \ = \ R_\lambda^\ast (\lambda-A)g \\
  &=& -P_{1,\lambda}(A-\lambda)^{-1} (\lambda-A)g \ =\ P_{1,\lambda} g=g  
\end{eqnarray*}
If we define $ \tilde{g}:=(\lambda-P_EA)g $, then it follows that $ \tilde{g}\in\overline{{\rm Ker}(B)\cap H^1} $.
Hence we have that $ R_\lambda^\ast\tilde{g}\in{\rm Dom}(A^\ast) $ and further $ (\lambda-A^\ast)R_\lambda^\ast\tilde{g}=\tilde{g} $.
As seen above, we have
$$    R_\lambda^\ast\tilde{g}=R_\lambda^\ast(\lambda-P_EA)g=g $$
which yields
$$   (\lambda-A^\ast)g=\tilde{g}=(\lambda-P_EA)g .$$
This finally proves that $ A^\ast g=P_EAg $. \hfill $ \Box $

We now want to investigate, whether the semigroups $ (T_t^{(c)})_{t\geq0} $ converge toward 
the semigroup $ (T_t^\ast)_{t\geq0} $. If $ (R^\ast_\lambda)_{\lambda>0} $ would be some true resolvent family, 
then this would be implied by Theorem 4.2 from Pazy (1983). 
The proof presented there however carries over to our situation. 
For convenience of the reader we include those arguments and adapt them to our particular situation. 

\begin{lemma}\label{Hilfslemma}
For every $ \varphi\in\overline{{\rm Ker}(B)\cap H^1} $, $ c\in\R $ and $ \lambda $ with $ {\rm Re}(\lambda)>0 $ we have
\begin{eqnarray*}
   R_\lambda^c\Big(T_t^{(c)}-T_t^\ast\Big)R_\lambda^\ast\varphi=\int_0^tT_{t-s}^{(c)}\Big(R_\lambda^\ast-R_\lambda^{(c)}\Big)T^\ast_s\varphi ds .
\end{eqnarray*}
\end{lemma}
{\bf Proof :}
For $ \varphi\in\overline{{\rm Ker}(B)\cap H^1} $  we note that the function 
$$ s\mapsto T_{t-s}^{(c)}R_\lambda^{(c)} T^\ast_sR_\lambda^\ast\varphi $$ 
is well defined. \\[0.5mm]
By differentiation we obtain
\begin{eqnarray*}
           \frac{d}{ds}T_{t-s}^{(c)}R_\lambda^{(c)} T^\ast_sR_\lambda^\ast\varphi
     &=&-T_{t-s}^{(c)}A_cR_\lambda^{(c)} T^\ast_sR_\lambda^\ast\varphi+T_{t-s}^{(c)}R_\lambda^{(c)} T^\ast_sA^\ast R_\lambda^\ast\varphi \\
       &=& T_{t-s}^{(c)}\Big(-A_cR_\lambda^{(c)}+R_\lambda^{(c)}A^\ast\Big) T^\ast_sR_\lambda^\ast\varphi \\
      &=& T_{t-s}^{(c)}\Big((\lambda-A_c)R_\lambda^{(c)}+R_\lambda^{(c)}(A^\ast-\lambda)\Big) T^\ast_s R_\lambda^\ast\varphi \\
         &=& T_{t-s}^{(c)}(\lambda-A_c)R_\lambda^{(c)} R_\lambda^\ast T^\ast_s\varphi
                                                         + T_{t-s}^{(c)}R_\lambda^{(c)}(A^\ast-\lambda) R_\lambda^\ast T^\ast_s\varphi \\
   &=& -T_{t-s}^{(c)}R_\lambda^\ast T^\ast_s\varphi + T_{t-s}^{(c)}R_\lambda^{(c)} T^\ast_s\varphi \\
        &=& -T_{t-s}^{(c)}\Big(R_\lambda^\ast - R_\lambda^{(c)} \Big) T^\ast_s\varphi .
\end{eqnarray*}
Now integration yields
\begin{eqnarray*}
       R_\lambda^{(c)}\Big( T^\ast_t-T_t^{(c)}\Big)R_\lambda^\ast\varphi
        &=&       R_\lambda^{(c)} T^\ast_tR_\lambda^\ast\varphi-T_t^{(c)}R_\lambda^{(c)}R_\lambda^\ast\varphi \\
        &=&    \Big[T_{t-s}^{(c)}R_\lambda^{(c)} T^\ast_sR_\lambda^\ast\varphi\Big]_0^t \\
        &=& -\int_0^tT_{t-s}^{(c)}\Big(R_\lambda^\ast - R_\lambda^{(c)} \Big) T^\ast_s\varphi ds.
\end{eqnarray*}
This finishes the proof of the lemma.
\hfill $ \Box $

\begin{theorem} \label{HalbgruppenKonvergenz}
For all   $ t\in[0,\infty[ $ and $ \psi\in\overline{{\rm Ker}(B)\cap H^1} $ one has $ T_t^{(c)}\psi $ converges toward $ T_t^\ast\psi $ in $ H $ as $ |c|\rightarrow\infty $. 
\end{theorem}
{\bf Proof :}
We know that for all $ \lambda>0 $ 
\begin{eqnarray} \label{Konv}
 R_\lambda^{(c)}\xi \longrightarrow R_\lambda^\ast\xi \ \ \ \mbox{for all $ \xi\in H $ as $ |c|\rightarrow\infty $.}
\end{eqnarray}
The triangle inequality yields for $  \varphi\in\overline{{\rm Ker}(B)\cap H^1} $ that
\begin{eqnarray*}
  \Big\|\Big(T_t^{(c)}-T_t^\ast\Big)R_\lambda^\ast\varphi\Big\| 
      &\leq&   \Big\|T_t^{(c)}\Big(R_\lambda^\ast-R_\lambda^{(c)}\Big)\varphi\Big\| +  \Big\|R_\lambda^{(c)}\Big(T_t^{(c)}-T_t^\ast\Big)\varphi\Big\| \\
      &&                                     +  \Big\|\Big(R_\lambda^{(c)}-R_\lambda^\ast\Big)T_t^\ast\varphi\Big\|   .
\end{eqnarray*}
For the first term we use the contraction property and (\ref{Konv}) to see that
$$  \Big\|T_t^{(c)}\Big(R_\lambda^\ast-R_\lambda^{(c)}\Big)\varphi\Big\|  \leq \Big\|\Big(R_\lambda^\ast-R_\lambda^{(c)}\Big)\varphi\Big\|\longrightarrow0 
                                                      \ \ \ \mbox{as $ |c|\rightarrow\infty $}.$$
For the third term we apply (\ref{Konv}) to the function $ \xi=T_t^\ast\varphi $ to obtain
$$  \Big\|\Big(R_\lambda^{(c)}-R_\lambda^\ast\Big)T_t^\ast\varphi\Big\| \longrightarrow0     \ \ \ \mbox{as $ |c|\rightarrow\infty $}.$$
In order to treat the second term we note that by Lemma \ref{Hilfslemma} for any choice of $ \psi\in H $ we have
\begin{eqnarray*}
     \Big\|R_\lambda^{(c)}\Big(T_t^{(c)}-T_t^\ast\Big) R_\lambda^\ast\psi\Big\|
    &\leq&  \int_0^t\Big\|T_{t-s}^{(c)}\Big(R_\lambda^\ast-R_\lambda^{(c)}\Big)T^\ast_s\psi\Big\| ds . 
\end{eqnarray*}
We note that one can use (\ref{Konv}) with $ \xi=T_s^\ast\psi $ to see that for every $ s\in[0,t] $ one has
\begin{eqnarray*}
   \Big\|T_{t-s}^{(c)}\Big(R_\lambda^\ast-R_\lambda^{(c)}\Big)T^\ast_s\psi\Big\| 
    \leq   \Big\|\Big(R_\lambda^\ast-R_\lambda^{(c)}\Big)T^\ast_s\psi\Big\| \longrightarrow0        \ \ \ \mbox{as $ |c|\rightarrow\infty $}.     
\end{eqnarray*}
Since moreover
\begin{eqnarray*}
   \Big\|T_{t-s}^{(c)}\Big(R_\lambda^\ast-R_\lambda^{(c)}\Big)T^\ast_s\psi\Big\| 
    &\leq&   \Big\|\Big(R_\lambda^\ast-R_\lambda^{(c)}\Big)T^\ast_s\psi\Big\|            \\
      &\leq&   \Big\|R_\lambda^\ast T^\ast_s\psi\Big\|+\Big\|R_\lambda^{(c)}T^\ast_s\psi\Big\|       \\
      &\leq&   \!\! \int_0^\infty \!\!\! e^{-\lambda r} \Big\|T_r^\ast T^\ast_s\psi\Big\|dr   \!
                                      +\! \int_0^\infty \! \!\! e^{-\lambda r} \Big\|T_r^{(c)}T^\ast_s\psi\Big\|dr \\
   &\leq&  2\int_0^\infty e^{-\lambda r} \|\psi\|dr   \\
  &=& \frac{2}{\lambda}\|\psi\| ,
\end{eqnarray*}
it follows from Lebesgue dominated convergence  theorem that 
\begin{eqnarray*}
     \Big\|R_\lambda^{(c)}\Big(T_t^{(c)}-T_t^\ast\Big) R_\lambda^\ast\psi\Big\| \longrightarrow0  \ \ \ \mbox{as $ |c|\rightarrow\infty $}.    
\end{eqnarray*}
Now, if we assume that $ \varphi\in{\rm Dom}(A^\ast) $, then it can be represented as $ \varphi=R_\lambda^\ast\psi $ with some suitable
$ \psi\in \overline{{\rm Ker}(B)\cap H^1} $. Then it follows that
$$   \Big\|R_\lambda^{(c)}\Big(T_t^{(c)}-T_t^\ast\Big)\varphi\Big\|
            = \Big\|R_\lambda^{(c)}\Big(T_t^{(c)}-T_t^\ast\Big) R_\lambda^\ast\psi\Big\| \longrightarrow0  \ \ \ \mbox{as $ |c|\rightarrow\infty $}.   $$
So far we have seen that for all $ \varphi\in {\rm Dom}(A^\ast) $ holds
$$      \Big\|\Big(T_t^{(c)}-T_t^\ast\Big)R_\lambda^\ast\varphi\Big\|   \longrightarrow0  \ \ \ \mbox{as $ |c|\rightarrow\infty $}. $$
Note also that every $ \vartheta\in {\rm Dom}((A^\ast)^2) $ can be represented as  $ \vartheta=R_\lambda^\ast\varphi $ with some suitable
$ \varphi\in{\rm Dom}(A^\ast) $. Therefore, we have proved that for all $  \vartheta\in {\rm Dom}((A^\ast)^2) $ holds
$$      \Big\|\Big(T_t^{(c)}-T_t^\ast\Big)\vartheta\Big\|   \longrightarrow0  \ \ \ \mbox{as $ |c|\rightarrow\infty $}. $$
According to Theorem 1.2.7 in Pazy (see Pazy p.6) we know that since $ T^\ast_t;t\geq0 $ is a $ C_0 $-semi group (i.e. : strongly continuous) its generator $ A^\ast $ satisfies that $ {\rm Dom}((A^\ast)^2) $ is dense in
$ \overline{{\rm Ker}(B)\cap H^1} $. It is now easy to use the boundedness of the semigroups to conclude that for any choice of $ \varphi\in \overline{{\rm Ker}(B)\cap H^1} $ one has
$$      \Big\|\Big(T_t^{(c)}-T_t^\ast\Big)\varphi\Big\|   \longrightarrow0  \ \ \ \mbox{as $ |c|\rightarrow\infty $} .$$
This finishes the proof.
\hfill $ \Box $

\section{The convergence of the diffusions}

In this section, we want to study the behavior of the diffusion processes, which are associated to the resolvent families $ R_\lambda^{(c)} $, as 
$ |c|\rightarrow\infty $.
According to Hsu (2001, p.24) there exists a probability space $ (\Omega,{\mathcal A},\mathbb{P}) $ a filtration $ {\mathcal F}_t,t\geq0 $ and 
$ {\mathcal F}_t $-adapted diffusion processes $ Y^{(c)}_t,t\geq0 $ with the property that for all $ f\in C^\infty(M) $ the process
$$  M_f^{(c)}(t) :=f(Y_t^{(c)})-f(Y^{(c)}_0)-\int_0^t A_cf(Y^{(c)}_s)ds $$
is a local $ {\mathcal F}_t $-martingale.
Moreover, it is proved in Hsu (2001 p.30) that the laws of those processes are uniquely determined.
Following Fukushima (1980, p.184) the regular Dirichlet form $ {\mathcal E}_* $ generates a symmetric Hunt-process $ Z_t;t\geq0 $ on the compact separable Hausdorff space $ X $ equipped with the Radon measure of full support $ \mu $.   \\[0.5mm]

For the following theorem we will need the following assumptions on the initial distribution :\\[0.5mm]
{\bf A1 :} For all $ c\in\R $ the initial distributions  $ \rho_0={\mathcal L}(Y^{(c)}_0) $ do not depend on $ c $ and that $ \rho_0 $ has a density $ g\in C(M) $ with respect to the measure $ \pi $ on $ M $.\\[0.5mm]
{\bf A2 :} The process  $ Z_0 $ is distributed according to the density $ \tilde{g}=G(P_I(g))  $ with respect to $ \mu $ on $ X $, where $ G $ is the isometry 
constructed in (\ref{GIsometry}) and $ P_I : L^2(M,\pi) \rightarrow  I^2(M) $ is the orthogonal projection, that we defined in the proof of Theorem \ref{Isometrie}.\\[0.5mm]

\begin{theorem} \label{CKonvergenz}
Under assumptions A1 and A2 we have for all $ f\in C(M)\cap \overline{{\rm Ker}(B)\cap H^1} $ that the finite dimensional distributions of processes $ f(Y^{(c)}_t);t\geq0 $ converges in distribution toward the finite dimensional distributions of the process 
$ Gf(Z_t),t\geq0 $, where $ G :\overline{{\rm Ker}(B)\cap H^1}\rightarrow L^2(X,\mu) $ was defined in (\ref{GIsometry}).
\end{theorem}
{\bf Proof :}
We first prove the convergence of the finite dimensional distributions.
According to Theorem \ref{HalbgruppenKonvergenz} we have for all $  f\in \overline{{\rm Ker}(B)\cap H^1}  $  that
$$  T_t^{(c)}f\longrightarrow T^\ast_tf \ \ \ \mbox{in $ H $ as $ |c|\rightarrow\infty $.} $$
If we denote as $ S_t;t\geq0 $ the semi-group on $ L^2(X,\mu) $ which is associated to the process $ Z_t;t\geq 0 $, then we have the following commutative diagram involving the two isometries $ F $ and $ G $, which where introduced in (\ref{FIsometry}) and (\ref{GIsometry}) :\\[1ex]
\begin{tikzpicture}
  \matrix (m) [matrix of math nodes,row sep=3em,column sep=4em,minimum width=2em] {
     \overline{{\rm Ker}(B)\cap H^1} &  \overline{{\rm Ker}(B)\cap H^1} \\
     L^2(X,\mu) & L^2(X,\mu) \\};
  \path[-stealth]
    (m-1-1) edge node [left] {$ G $} (m-2-1)
            edge node [above] {$ T_t^\ast $} (m-1-2)
    (m-2-1.east|-m-2-2) edge node [above] {$ S_t $} (m-2-2)
    (m-2-2) edge node [left] {$ F $} (m-1-2);
\end{tikzpicture}

For $ f\in C^2(M)\cap\overline{{\rm Ker}(B)\cap H^1}\subset C^\sim(M)\cap\overline{{\rm Ker}(B)\cap H^1} $ there is a $ \tilde{f}\in C(X) $ such that $ f=\tilde{f}\circ p_X $ and thus $ G(f)=\tilde{f} $.  
In that situation we also have that 
$$ \exp(i\xi f)=\exp(i\xi\tilde{f})\circ p_X\in C^\sim(M)\cap\overline{{\rm Ker}(B)\cap H^1}  $$ 
and thus it follows that
$$  G(\exp(i\xi f))=\exp(i\xi\tilde{f})=\exp(i\xi G(f)) .$$
If we use the fact that $ F $ is the inverse of $ G $ and that both are $ L^2 $-isometries (see Theorem \ref{Isometrie}), this implies for $ \xi\in\mathbb{R} $ that
\begin{eqnarray*}
   \lim_{|c|\rightarrow\infty}\mathbb{E}\left[\exp(i\xi f(Y_t^{(c)}))\right]&=&\lim_{|c|\rightarrow\infty}\int_M gT_t^{(c)}\exp(i\xi f)d\pi\\
       &=& \int_M gT_t^\ast\exp(i\xi f)d\pi    \\ 
       &=&  \int_M P_I(g)T_t^\ast\exp(i\xi f)d\pi  \ \ \ \mbox{(since $ T_t^\ast\exp(i\xi f) \in I^2(M) $)}\\
        & = &  \int_M P_I(g) F\left( S_t\left( G\left(\exp(i\xi f)\right)\right)\right)d\pi \\
            &=&\int_X \tilde{g}S_t\!\left(\exp(i\xi Gf)\right)d\mu \\
            & = &   \mathbb{E}\Big[\exp(i\xi G(f)(Z_t)\Big]  .
\end{eqnarray*}
This proves that the one dimensional marginal distribution converge. 
In order to treat the higher dimensional case, we use the Markov property to see that for $ \xi_1,\xi_2\in\mathbb{R} $ we have
\begin{eqnarray*}
   && \mathbb{E}\Big[\exp\left(i\xi_1f(Y_{t_1}^{(c)})+i\xi_2f(Y_{t_2}^{(c)})\right)\Big] \\ 
   &=&    \mathbb{E}\left[\exp\left(i\xi_1f(Y_{t_1}^{(c)})\right)\mathbb{E}_{Y_{t_1}^{(c)}}\left[\exp\left(i\xi_2f(Y_{t_2-t_1}^{(c)})\right)\right]\right]  \\
   &=&  \int_M gT_{t_1}^{(c)}\left(e^{i\xi_1f}T_{t_2-t_1}^{(c)}e^{i\xi_2f}\right)d\pi .
\end{eqnarray*}
Further, since by Proposition \ref{Homeomorphism} the map $ G $ is an algebra homeomorphism, it follows that
\begin{eqnarray*}
  && \mathbb{E}\Big[\exp\big(i\xi_1G(f)(Z_{t_1})+i\xi_2G(f)(Z_{t_2})\big)\Big] \\
  &=&   \int_X \tilde{g}\ S_{t_1}\!\left(e^{i\xi_1G(f)}S_{t_2-t_1} e^{i\xi_2G(f)}\right)d\mu  \\
  &=&   \int_X G(P_Ig)\ S_{t_1}\!\left(e^{i\xi_1G(f)}G\big(T^\ast_{t_2-t_1} e^{i\xi_2f}\big)\right)d\mu  \\
  &=&   \int_X G(P_Ig)\ S_{t_1}\!\Big(G\big(e^{i\xi_1f}T^\ast_{t_2-t_1} e^{i\xi_2f}\big)\Big)d\mu  \\
  &=&   \int_X G(P_Ig)\ G\Big(T^\ast_{t_1}\!\big(e^{i\xi_1f}T^\ast_{t_2-t_1} e^{i\xi_2f}\big)\Big)d\mu  \\
  &=&   \int_X P_Ig\ \Big(T^\ast_{t_1}\!\big(e^{i\xi_1f}T^\ast_{t_2-t_1} e^{i\xi_2f}\big)\Big)d\pi  \\
   &=&  \int_M gT_{t_1}^\ast\left(e^{i\xi_1f}T_{t_2-t_1}^\ast e^{i\xi_2f}\right)d\pi .
\end{eqnarray*}
Furthermore, for a given $ \epsilon>0 $ as $ |c| $ becomes sufficiently large, we have
\begin{eqnarray*}
 &&  \Bigg|\int_M gT_{t_1}^\ast\Big(e^{i\xi_1f}T^\ast_{t_2-t_1}e^{i\xi_2f}\Big)d\pi
                 -\int_M gT_{t_1}^{(c)}\Big(e^{i\xi_1f}T^{(c)}_{t_2-t_1}e^{i\xi_2f}\Big)d\pi\Bigg| \\
 &\leq& \Bigg|\int_M gT_{t_1}^\ast\Big(e^{i\xi_1f}T^\ast_{t_2-t_1}e^{i\xi_2f}\Big)d\pi
                 -\int_M gT_{t_1}^{(c)}\Big(e^{i\xi_1f}T^\ast_{t_2-t_1}e^{i\xi_2f}\Big)d\pi\Bigg| \\
&&+\Bigg|\int_M gT_{t_1}^{(c)}\Big(e^{i\xi_1f}T^\ast_{t_2-t_1}e^{i\xi_2f}\Big)d\pi
                 -\int_M gT_{t_1}^{(c)}\Big(e^{i\xi_1f}T^{(c)}_{t_2-t_1}e^{i\xi_2f}\Big)d\pi\Bigg| \\
&\leq& \epsilon+\Bigg|\int_M gT_{t_1}^{(c)}\Big(e^{i\xi_1f}T^\ast_{t_2-t_1}e^{i\xi_2f}
                                        -   e^{i\xi_1f}T^{(c)}_{t_2-t_1}e^{i\xi_2f}\Big)d\pi\Bigg| .
\end{eqnarray*}
By Cauchy-Schwarz ineqality we have that the second term in the last line is bounded by the square root of the following expression
\begin{eqnarray*}
  && \|g\|\Big\|T_{t_1}^{(c)}\Big(e^{i\xi_1f}T^\ast_{t_2-t_1}e^{i\xi_2f} -   e^{i\xi_1f}T^{(c)}_{t_2-t_1}e^{i\xi_2f}\Big)\Big\| \\
  &\leq& \|g\|\Big\|e^{i\xi_1f}T^\ast_{t_2-t_1}e^{i\xi_2f} -   e^{i\xi_1f}T^{(c)}_{t_2-t_1}e^{i\xi_2f}\Big\| \\
  &\leq& \|g\|\Big\|T^\ast_{t_2-t_1}e^{i\xi_2f} -  T^{(c)}_{t_2-t_1}e^{i\xi_2f}\Big\| \longrightarrow0 \ \ \mbox{as $ |c|\rightarrow\infty $.}
\end{eqnarray*}
This proves that as $ |c|\rightarrow\infty $ the bivariate characteristic distribution
$$   \phi_{f(Y_{t_1}^{(c)}),f(Y_{t_2}^{(c)})}(\xi_1,\xi_2)=\mathbb{E}\left[\exp\left(i\xi_1G(f)(Y_{t_1}^{(c)})+i\xi_2G(f)(Y_{t_2}^{(c)})\right)\right]    $$
converge toward the characteristic function
$$   \phi_{G(f)(Z_{t_1}),G(f)(Z_{t_2})}(\xi_1,\xi_2) =\mathbb{E}\left[\exp\left(i\xi_1G(f)(Z_{t_1})+i\xi_2G(f)(Z_{t_2})\right)\right] .$$
Iteration of those arguments proves the convergence of the all finite dimensional distributions.  
\hfill $ \Box $

\begin{corollary} \label{CCorollary}
Under assumptions A1 and A2 the finite dimensional distributions of the equivalence class valued  processes $ [Y^{(c)}_t];t\geq0 $ converge weakly toward the finite dimensional distributions of the process $ Z_t;t\geq0 $.
\end{corollary}
{\bf Proof :} This is a direct consequence of Theorem \ref{CKonvergenz}.\hfill $ \Box $

\begin{theorem} \label{C2Konvergenz}
Under assumptions A1 and A2 we have that for all $ f\in C^2(M)\cap \overline{{\rm Ker}(B)\cap H^1} $ that the family of processes $ f(Y^{(c)}_t);t\geq0 $ converges in distribution toward the process 
$ Gf(Z_t),t\geq0 $, where $ G :\overline{{\rm Ker}(B)\cap H^1}\rightarrow L^2(X,\mu) $ was defined in (\ref{GIsometry}).\footnote{Does this prove that $ Z $ has continuous trajectories ?}
\end{theorem}
{\bf Proof :}
It is now sufficient to prove tightness of the family of processes $ f(Y_t^{(c)}),t\geq0 $.
According to Billingsley (1968 p.55) the sequence of processes $ f(Y^{(c)}_t);t\in[0,1] $ is tight in $ C_\R([0,1]) $ when the two following conditions are satisfied :\\[0.5mm]
i) the family of random variables $ f(Y^{(c)}_0), c\in\R $ is tight in $ \R $;\\[0.5mm]
ii) for all $ \eta>0 $ and $ \epsilon>0 $ there exists a $ \delta>0 $ and a $ c_0>0 $ with the property 
$$   \mathbb{P}\Bigg( \sup_{|s-t|<\delta}\big|f(Y_s^{(c)})-f(Y_t^{(c)})\big|\geq\epsilon \Bigg)\leq \eta \ \ \ \mbox{for all $ c $ with $ |c|>c_0 $.}  $$
The first condition holds, since the $ C^2 $-function $ f $ is bounded on the compact manifold $ M $.\\[0.5mm] 
For the second condition, we note that
\begin{eqnarray*}
 &&  \mathbb{P}\Bigg( \sup_{|s-t|<\delta}\big|f(Y_s^{(c)})-f(Y_t^{(c)})\big|\geq\epsilon \Bigg) \\
 &\leq&  \mathbb{P}\Bigg( \sup_{|s-t|<\delta}\big|M_f^{(c)}(t)-M_f^{(c)}(s)\big|+\sup_{|s-t|<\delta}\int_s^t \big|A_cf(Y^{(c)}_r)\big|dr\geq\epsilon \Bigg) \\
 &\leq&  \mathbb{P}\Bigg( \sup_{|s-t|<\delta}\big|M_f^{(c)}(t)-M_f^{(c)}(s)\big|\geq\frac{\epsilon}{2}\Bigg)
                                            +\mathbb{P}\Bigg(\sup_{|s-t|<\delta}\int_s^t \big|A_cf(Y^{(c)}_r)\big|dr\geq\frac{\epsilon}{2} \Bigg) .
\end{eqnarray*}
By compactness of $ M $ for $ f\in  C^2(M)\cap \overline{{\rm Ker}(B)\cap H^1} $ there exists a $ C_f>0 $ such that 
$$  \big| A_cf\big|=\big| A_0f\big|=\big|\Delta f -\nabla U\cdot\nabla f\big|\leq C_f .$$
Then the second term yields
\begin{eqnarray*}
\mathbb{P}\Bigg(\sup_{|s-t|<\delta}\int_s^t \big|A_cf(Y^{(c)}_r)\big|dr\geq\frac{\epsilon}{2} \Bigg)
   \leq  \mathbb{P}\Bigg(\sup_{|s-t|<\delta}C_f|t-s|\geq\frac{\epsilon}{2} \Bigg) =  \mathbb{P}\Big(C_f\delta\geq\frac{\epsilon}{2} \Big) ,
\end{eqnarray*}
which becomes zero for $ \delta<\frac{\epsilon}{2C_f} $.\\[0.5mm]
In order to obtain some bound on the first term, we need some more insight into the construction of the process $ Y^{(c)} $. 
The definition of the process $ Y^{(c)} $ uses the fact that the $ d $-dimensional Riemannian manifold $ (M,g) $ can be embedded into 
some euclidean space $ \R^N $ by virtue of Whitney's theorem. In local coordinates the differential operator $ A_c $ takes the form :
$$   A_cf=\frac{1}{2}\sum_{i,j=1}^da_{ij}\partial_{x_i}\partial_{x_j}f+\sum_{i=1}^db_i\partial_{x_i}f   ,$$
where the matrix valued function $ (a_{ij})_{1\leq i,j\leq d} $ is symmetric. 
The coefficient functions $ a_{ij}, 1\leq i,j\leq d $ and $ b_i, 1\leq i\leq d $ can be extended to some functions  
$ \tilde{a}_{ij}, 1\leq i,j\leq d $ and $ \tilde{b}_i, 1\leq i\leq d $ defined on some neighborhood of $ M $ in $ \R^N $. The resulting
differential operator 
$$   \tilde{A}_cf=\frac{1}{2}\sum_{i,j=1}^d\tilde{a}_{ij}\partial_{x_i}\partial_{x_j}f+\sum_{i=1}^d\tilde{b}_i\partial_{x_i}f    $$
satisfies $ \tilde{A}_c\tilde{f}=A_cf $ whenever $ \tilde{f} $ is a twice continously differentiable extension of a function $ f\in C^2(M) $ to $ \R^N $. 
We note $ \tilde{a} $ the matrix valued function $ (a_{ij})_{1\leq i,j\leq d} $ and $ \tilde{b} $ the vector valued function $ (\tilde{b}_i)_{1\leq i\leq d} $.
We note that there exists a unique symmetric positive definite matrix  $ \tilde{\sigma} $ satisfying $ \tilde{\sigma}\tilde{\sigma}^t=\tilde{a} $. 
The operator $ \tilde{A}_c $ is the generator of the $ \R^N $-valued diffusion process $ Y^{(c)} $ which is the solution of the following
stochastic differential equation
$$   dY^{(c)}_t=\tilde{b}(Y^{(c)}_t)dt+\tilde{\sigma}(Y^{(c)}_t) dW_t   $$
where $ W_t,t\geq0 $ is standard Brownian motion on $ \R^N $ on a suitable filtered probability space  $ (\Omega,{\mathcal A},\mathbb{P}, ({\mathcal F}_t)_{t\geq0}) $  . 
It is proved in Hsu (2001 p.28) that the diffusion $ Y^{(c)} $ stays in $ M $ if $ Y^{(c)}_0\in M $.
In order to obtain a representation of the local martingale $ M_f^{(c)} $ in terms of stochastic integrals we use It\^{o}-formula to see that 
for some $ \tilde{f}\in C^\infty(M) $ we have 
\begin{eqnarray*}
     \tilde{f}(Y^{(c)}_t) &=&   \int_0^t\sum_{i=1}^N\partial_{x_i}\tilde{f}(Y^{(c)}_s)\tilde{\sigma}_{ij}(Y^{(c)}_s)dW^j_s+R_t ,
\end{eqnarray*}
where $ R_t,t\geq0 $ has bounded variation. Thus we have that 
\begin{eqnarray*}
   M_f^{(c)}(t)=  \int_0^t\sum_{i,j=1}^N\partial_{x_i}\tilde{f}(Y^{(c)}_s)\tilde{\sigma}_{ij}(Y^{(c)}_s)dW^j_s .
\end{eqnarray*}
By Nash's embedding theorem we can assume without loss of generality that the embedding used to construct $ Y^{(c)} $ is a Riemannian isometry
mapping $ (M,g) $ into $ \R^N $ with its usual Riemannian structure. Thus the functions $ \tilde{f} $, $ \tilde{\sigma}_{ij},1\leq i,j\leq N $ can be assumed
to have compact support. Thus those functions can be assumed to be bounded by some constant. 
Following some argument that we found in Bensoussan et al. (1978 p.403), we note that for 
$$ Q :=  \sup_{ 0\leq s,t\leq 1, |s-t|<\delta}\big|M_f^{(c)}(t)-M_f^{(c)}(s)\big|  $$
we have that
\begin{eqnarray*}
 Q  &=&  \sup_{ 0\leq s,t\leq 1, |s-t|<\delta}\Bigg|   \int_s^t\sum_{i,j=1}^N\partial_{x_i}\tilde{f}(Y^{(c)}_s)\tilde{\sigma}_{ij}(Y^{(c)}_s)dW^j_s   \Bigg| \\
 &\leq&  2 \sup_{k\in\{0,...,[\frac{1}{\delta}]\}}\sup_{k\delta\leq t\leq (k+1)\delta}
                               \Bigg|   \int_{k\delta}^t\sum_{i,j=1}^N\partial_{x_i}\tilde{f}(Y^{(c)}_s)\tilde{\sigma}_{ij}(Y^{(c)}_s)dW^j_s   \Bigg| .
\end{eqnarray*}
It follows that
\begin{eqnarray*}
  \mathbb{P}\Bigg( Q\geq\frac{\epsilon}{2}\Bigg)   \leq \sum_{k\in\{0,...,[\frac{1}{\delta}]\}}    \! \! \! \! \mathbb{P}\Bigg(\sup_{k\delta\leq t\leq (k+1)\delta}
                      \Bigg|   \int_{k\delta}^t\sum_{i,j=1}^N\partial_{x_i}\tilde{f}(Y^{(c)}_s)\tilde{\sigma}_{ij}(Y^{(c)}_s)dW^j_s   \Bigg| \geq \frac{\epsilon}{4}\Bigg) .
\end{eqnarray*}
Applying first Chebyshev inequality and then Doob inequality (see Karatzas and Shreve (1987 p.14)) yields
\begin{eqnarray*}
  &&  \mathbb{P}\Bigg(   \sup_{ 0\leq s,t\leq 1, |s-t|<\delta}\big|M_f^{(c)}(t)-M_f^{(c)}(s)\big|   \geq\frac{\epsilon}{2}\Bigg) \\
    &\leq&  \frac{4^4}{\epsilon^4}\sum_{k\in\{0,...,[\frac{1}{\delta}]\}}
         \mathbb{E}\Bigg[ \sup_{k\delta\leq t\leq (k+1)\delta}  \Bigg|   
                      \int_{k\delta}^t\sum_{i,j=1}^N\partial_{x_i}\tilde{f}(Y^{(c)}_s)\tilde{\sigma}_{ij}(Y^{(c)}_s)dW^j_s   \Bigg|^4\Bigg] \\
  &\leq&  \frac{4^9}{(3\epsilon)^4}\sum_{k\in\{0,...,[\frac{1}{\delta}]\}}
         \mathbb{E}\Bigg[ 
                        \Bigg|  \int_{k\delta}^{(k+1)\delta}\sum_{i,j=1}^N\partial_{x_i}\tilde{f}(Y^{(c)}_s)\tilde{\sigma}_{ij}(Y^{(c)}_s)dW^j_s   \Bigg|^4\Bigg] .
\end{eqnarray*}
It now follows from Burkholder-Davis-Gundy inequality that there exists a constant $ C>0 $ such that for all $ k $ one has
\begin{eqnarray*}
 &&  \mathbb{E}\Bigg[ 
                        \Bigg|  \int_{k\delta}^{(k+1)\delta}\sum_{i,j=1}^N\partial_{x_i}\tilde{f}(Y^{(c)}_s)\tilde{\sigma}_{ij}(Y^{(c)}_s)dW^j_s   \Bigg|^4\Bigg] \\
     &\leq& C \mathbb{E}\Bigg[ 
    \Bigg(  \int_{k\delta}^{(k+1)\delta}\sum_{i=1}^N\Big(\sum_{j=1}^d\partial_{x_i}\tilde{f}(Y^{(c)}_s)\tilde{\sigma}_{ij}(Y^{(c)}_s)\Big)^2ds   \Bigg)^2\Bigg] \\
                          &\leq& CE_f\frac{\delta^2}{\epsilon^4} ,
\end{eqnarray*}
where the constant $ E_f $ only depends on the bounds for the coefficients $ \tilde{\sigma}_{ij} $ and the partial derivatives of $ \tilde{f} $. 
Finally this shows that 
\begin{eqnarray*}
   \mathbb{P}\Bigg(   \sup_{ 0\leq s,t\leq 1, |s-t|<\delta}\big|M_f^{(c)}(t)-M_f^{(c)}(s)\big|   \geq\frac{\epsilon}{2}\Bigg) 
 \ \leq \ CE_f\frac{\delta^2}{\epsilon^4}\sum_{k\in\{0,...,[\frac{1}{\delta}]\}}  1 
  \ \sim \ CE_f\frac{\delta}{\epsilon^4}  .
\end{eqnarray*}
This finishes the proof of the tightness since $ \delta $ can be choosen arbitrarily small.
\hfill $ \Box $

\begin{corollary}
Suppose that the assumptions A1 and A2 hold. 
If for all $ x\in M $ the functions $ y\mapsto d_X([x],[y]) $ are in $ C^2(M) $, then the equivalence class valued  processes $ [Y^{(c)}_t];t\geq0 $ converges in distribution toward the process $ Z_t;t\geq0 $.
\end{corollary}
{\bf Proof :} By Corollary \ref{CCorollary} it is sufficient to prove tightness for the family of process $ [Y^{(c)}_t];t\geq0 $. The tightness criterion from Billingsley in this context becomes :\\[0.5mm]
i) the family of random variables $ [Y^{(c)}_0], c\in\R $ is tight in $ X $;\\[0.5mm]
ii) for all $ \eta>0 $ and $ \epsilon>0 $ there exists a $ \delta>0 $ and a $ c_0>0 $ with the property 
$$   \mathbb{P}\Bigg( \sup_{|s-t|<\delta}d_X\big([Y_s^{(c)}],[Y_t^{(c)}]\big)\Big|\geq\epsilon \Bigg)\leq \eta \ \ \ \mbox{for all $ c $ with $ |c|>c_0 $.}  $$
The first condition follows from the compactness of the metric space $ X $. For the second condition, we note that by triangle inequality
\begin{eqnarray*}
  d_X\big([Y_s^{(c)}],[Y_t^{(c)}]\big)  &\leq& \sum_{i=0}^n d_X\big([Y_{r_i}^{(c)}],[Y_{r_{i+1}}^{(c)}]\big) ,
\end{eqnarray*}
where $ s=r_0<r_1<...<r_n=t $ is an arbitrary partition of $ [s,t] $. Refining this partition, we obtain as $ \max_i|r_{i+1}-r_i|\rightarrow0 $ the following
inequality
\begin{eqnarray*}
  d_X\big([Y_s^{(c)}],[Y_t^{(c)}]\big)  &\leq&  M_d^{(c)}(t)-M_d^{(c)}(s)+\int_s^t \big|A_cG_{[Y_s^{(c)}]}([Y_s^{(c)}])\big|ds ,
\end{eqnarray*}
where $ M_d^{(c)}(t) $ is a $ \mathcal{F}_t $-martingale and 
$$ G_x:M\rightarrow\mathbb{R};y\mapsto=d_X([x],[y])  . $$
From this inequality the arguments described in the proof of Theorem \ref{C2Konvergenz} can be used again. 
\hfill $ \Box $

\section{Examples} \label{Beispiele}

In this section we want to illustrate the main theorem through some examples. We begin with the two extremal situations.\\[1ex]
{\bf Example 1 :}
If the vectorfield $ b $ is the zero vectorfield, then $ \Phi_t(x)=x $ for all $ t\in\R $. The equivalence classes introduced in the previous section then only contain one element; i.e. : $ [x]=\{x\} $ for all $ x\in M $. From this results $ X=M $ and the generator of the limit diffusion is $ A_0 $.\\[1ex]
{\bf Example 2 :} 
Let $ b $ be a $ C^1 $-vectorfield such that the flow $ (\Phi_t)_{t\in\R} $ generated by the equation $ \dot{x}=b(x) $ has the property that for some 
$ x\in M $ one has $ \overline{\{\Phi_t(x);t\in\R\}} = M $. This property holds in particular, when the flow $ (\Phi_t)_{t\in\R} $ is transitive; i.e. : 
 $$ \exists x\in M : \ \ \ \left\{y :\liminf_{t\rightarrow\infty} d(y,\Phi_t(x))=0\right\} = M .$$
We note that every ergodic flow is also transitive (see Man\'{e} p.104).\\[0.5mm]
Since there exists a trajectory, which is dense, it follow that $ x\sim y $ for all $ x,y\in M $. Thus we have $ [x]=M $ for all $ x\in M $ and
 it follows that the state space $ X $ consists of only one point and the Hunt-process $ Y^\ast $ is trivial. \\[1ex]
 
{\bf Example 3 :} \label{ergodicperiodic}
The three dimensional torus $ T^3 $ can be represented by the cube $ [0,2\pi]^3 $ with the usual identifications
$  (x_1,x_2,x_3)=(y_1,y_2,y_3) $ iff $ x_i=0 $ and $ y_i=2\pi $ for at least one $ i\in\{1,2,3\} $.
On $ T^3 $ we define a divergence-free vectorfield through 
$$  b(x)=\big(0,\sin x_1,\cos x_1\big) .$$
It is not difficult to see that for all $ x_1\in[0,2\pi] $ the two dimensional tori 
$$  T^2_{x_1}=\Big\{(x_1,y_2,y_3);y_2,y_3\in[0,2\pi]\Big\} $$
are invariant under the flow $ (\Phi_t)_{t\in\R} $ generated by the vectorfield $ b $. 
It is well known, that the flow restricted to those two dimensional tori is ergodic iff $ \tan x_1 $ is irrational. Those points form a dense set in $ [0,2\pi] $.
In that case the flow is also transient and there exists some dense trajectory. As a consequence, we see that if $ \tan x_1 $ is irrational, then any  two elements $ x $, $ y $ from $ T^2_{x_1} $  are equivalent.
When $ \tan x_1 $ is rational and $ x $, $ y $ are two points from $ T^2_{x_1} $, then there exist two points $ u $ and $ v $ such that $  u_1= v_1 $, $ \tan u_1 $ is irrational, $ d(u,x)<\epsilon/3 $ and $ d(y,v)<\epsilon/3 $. Obviously, one can follow the flow-line through $ u $ until one comes $ \epsilon/3 $-close to $ v $. This reasoning shows, that $ x $ and $ y $ are equivalent. It follows that $ M=[0,2\pi] $ where $ 0 $ and $ 2\pi $ are identified.
The limiting diffusion is a Brownian motion on the circle.\\[1ex]

{\bf Example 4 :} Let $ S^2 $ be the two dimensional sphere. We obtain a parametrization of $ S^2 $ in $ \R^3 $ through 
$$ (\theta_1,\theta_2)\mapsto\Big(\cos\theta_1\sin\theta_2,\cos\theta_1\cos\theta_2,\sin\theta_1\Big) $$
Given a $ C^1 $-function $ \rho :]0,\pi]\rightarrow\R  $ the incompressible flow 
\[ \phi_t\left(\begin{array}{c}\cos\theta_1\sin\theta_2 \\
                                   \cos\theta_1\cos\theta_2 \\      
        \sin\theta_1           \end{array}\right)
       :=\left(\begin{array}{c}\cos\theta_1\sin(\theta_2+t\rho(\theta_1)) \\
                                   \cos\theta_1\cos(\theta_2+t\rho(\theta_1)) \\      
        \sin\theta_1           \end{array}\right) . \]
After differentiation with respect to $ t $ we see that this flow is generated by the following vectorfield on $ S^2 $
\[  b\left(\begin{array}{c}\cos\theta_1\sin\theta_2 \\
                                  \cos\theta_1\cos\theta_2 \\      
        \sin\theta_1           \end{array}\right)
  :=\rho(\theta_1)\left(\begin{array}{c}\cos\theta_1\cos\theta_2 \\
                                   -\cos\theta_1\sin\theta_2 \\      
        0    \end{array}\right) .\]
The flow $ (\Phi_t)_{t\in\R} $ generated from this vectorfield leaves invariant the small circles defined through 
$$  S^1_{\theta_1} :=\Big\{\Big(\cos\theta_1\sin\theta_2,\cos\theta_1\cos\theta_2,\sin\theta_1\Big);\theta_2\in[0,2\pi)\Big\} . $$ 
Since those small circles are parallel to each other it follows that any two elements from two different small circles can not be in the same 
equivalence class with respect to the relation defined in section \ref{Sektion5}. 
The flow generated by the vectorfield $ b $ is essentially a rotation along those small circles.
The geometry of the metric space $ X $ is determined through the set where the function $ \rho $ vanishes as follows :\\[0.5mm]
If $ \rho(\theta_1)\neq0 $ then all of $ S^1_{\theta_1} $ forms one equivalence class (since one can move along the circle at zero cost);\\[0.5mm]
If $ \rho(\theta_1)=0 $ and there exists a sequence of values $ \theta_1^{(n)} $ with $ \rho(\theta_1^{(n)})\neq0 $ converging toward $ \theta_1 $, then all of $ S^1_{\theta_1} $ forms one equivalence class (since one can reach a zero cost trajectory at arbitrary low price);\\[0.5mm]
If $ \theta_1 $ is in an open set, where the function $ \rho $ is zero, then any two elements from $ S^1_{\theta_1} $ are not equivalent.\\[0.5mm]
If $ \theta_1 $ is not in an open set, where the function $ \rho $ is zero, then all the elements from the small circle $ \theta_1 $ are 
in the same equivalence class.\\[0.5mm]
The resulting space $ X $ is thus obtained from the manifold $ M $ through collapsing of small circles where the function $ \rho $ is non-zero.\\[1ex]

{\bf Example 5 :} For $ 0<a_0<a_1<a_2 $ the ellipsoid is defined through
$$  E :=\Bigg\{\frac{x_0^2}{a_0}+\frac{x_1^2}{a_1}+\frac{x_2^2}{a_2}=1\Bigg\}. $$
The sphere boundle $ SE $ carries a natural Riemannian metric. The associated Laplace operator $ \Delta_{SE} $ 
generates Brownian motion on this resulting four dimensional compact Riemannian manifold. The geodesic flow  on $ SE $ is an integrable 
dynamical system which preserves the Liouville measure on $ SE $. It is generated by a divergence free vectorfield that we will denote by 
$ b_g $. In Franke, Hwang, Pai, Sheu (2010) one can find further information on the diffusions generated by the elliptic operators
$$     A_c :=\Delta_{SE}+cb_g\cdot\nabla .$$

\bigskip

{\bf Example 6 :} We present an example that shows that a continuous $ H^1 $ function can be constant along the flow lines without being constant along the equivalence classes with respect to the relation $ \sim $ introduced for the construction of the state space $ X $. Let $ C $ be the Cantor set in $ [0,1] $. In the usual construction of the Cantor set one starts with the interval $ C_0=[0,1] $ from which one removes the open interval $ G_{1,1}=]1/3,2/3[ $. Let $ x_{1,1} $ be the center point of this interval. From the resulting set $ C_1=[0,1]\backslash G_{1,1} $ one removes the two open intervals $ G_{2,1}=]1/9,2/9[ $ and $ G_{2,2}=]7/9,8/9[ $ one obtains the set $ C_2=C_1\backslash(G_{2,1}\cup G_{2,2}) $. Let $ x_{2,1}=1/6 $ and $ x_{2,2}=5/6 $ be the mid points of the two intervals $ G_{2,1} $ resp. $ G_{2,2} $. We can continue to obtain
$$ C_3=C_2\backslash(G_{3,1}\cup G_{3,2}\cup G_{3,3}\cup G_{3,4}) $$ 
with center points $ x_{3,1}=1/18 $, $ x_{3,2}=5/18 $, $ x_{3,3}=13/18 $ and $ x_{3,4}=17/18 $ for the intervals $ G_{3,1} $, $ G_{3,2} $, $ G_{3,3} $ and $ G_{3,4}$. In the $ n $-th step of the construction, one has
$$ C_n=C_{n-1}\backslash(G_{n,1}\cup ...\cup G_{n,2^{n-1}}) $$ and
center points $ x_{n,1},...,x_{n,2^{n-1}} $ for the intervals $ G_{n,1},..., G_{n,2^{n-1}} $. The Cantor set is then obtained as the intersection $$ C:=\bigcap_{n\in\mathbb{N}}C_n .$$
Note that the distances between one of those points $ x_{n,k} $ and the Cantor set $ C $ equals $ r_n:=1/(2\cdot3^{n}) $. 
Further, it is well known, that the Lebesgue measure of the sets $ C_n $ converge to zero as $ n\rightarrow\infty $.
The Cantor function $ f:[0,1]\rightarrow[0,1] $ is continuous, takes constant values on all gaps $ C_{n,k} $ and satisfies $ f(0)=0 $ and $ f(1)=1 $ . 
It can be obtained as a uniform limit of a sequence of continuous functions $ f_n $. 
The construction of the sequence $ f_n;n\in\mathbb{N} $ is done through induction. 
One starts with $ f_0(x)=x $ for all $ x\in[0,1] $. For $ n\geq1 $ functions $ f_n $ are defined inductively as follows
\[  
              f_{n+1}(x)=\left\{\begin{array}{cl}  \frac{1}{2}f_n(3x) & \mbox{for $ 0\leq x\leq \frac{1}{3} $;} \\
                                                               \frac{1}{2} & \mbox{for $ \frac{1}{3}< x\leq \frac{2}{3} $;} \\
                                                             \frac{1}{2}+\frac{1}{2}f_n(3x-2)       & \mbox{for $ \frac{2}{3}< x\leq 1 $.} \end{array}\right.
\]
By construction, the function $ f_n $ is continuous, is constant on the gaps $ C_{m,k} $ for all $ m\leq n $ and satisfies $ f_n(0)=0 $ and $ f_n(1)=1 $. The sequence $ f_n;n\in\mathbb{N} $ converges uniformly toward a function $ f $ as $ n\rightarrow\infty $. It follows that this limit function $ f $ is continuous, constant on the gaps of the Cantor set and satisfies $ f(0)=0 $ and $ f(1)=1 $. 
Let $ T^2 $ be the two dimensional torus, which can be constructed from the square $ [-1,1]\times[-1,1] $ with the following relations $ (-1,y)\sim(1,y) $ and $ (x,-1)\sim(x,1) $. On the torus, we can define the family of points $ z_{n,k}=(x_{n,k},0) $ for $ n\geq1 $ and $ 1\leq k\leq 2^{n-1} $. The above consideration then shows that the family of open balls $ B_{r_n}(z_{n,k}) $; $ n\geq1 $ and $ 1\leq k\leq 2^{n-1} $ are pairwise disjoint and if one removes their intersections from the geodesic circle $ \{(x,0);x\in[-1,1]\} $ in $ T^2 $ one obtains the Cantor set on the segment $  \{(x,0);x\in[0,1]\}  $. Let $ \rho:\mathbb{R}\rightarrow\mathbb{R} $ be a smooth function which vanishes outside $ [1/2,1] $ and which is positive for all $ x\in ]1/2,1[ $. For a given $ z_{n,k}=(x_{n,k},0) $ we can construct a vectorfield $ b_{n,k} $ with support in $ B_{r_n}(z_{n,k}) $ as follows
\[
  b_{n,k}(x,y)=\rho\Big(r_n^{-1}\sqrt{(x-x_{n,k})^2+y^2}\Big)\left(\begin{array}{c} y  \\ x_{n,k}-x
                  \end{array}\right).
\]
This vector field generates circular flow lines with center $ z_{n,k} $ for any radius $ r\in]r_n/2,r_n[ $.

\begin{tikzpicture}
\draw (0,0) -- (10,0);
\draw (10,0) -- (10,10);
\draw (10,10) -- (0,10);
\draw (0,10) -- (0,0);
\draw (0+5,5) -- (5/81+5,5);
\draw (10/81+5,5) -- (5/27+5,5);
\draw (10/27+5,5) arc (0:350:5/54);
\draw (10/27+5,5) -- (35/81+5,5);
\draw (40/81+5,5) -- (5/9+5,5);
\draw[->] (10/9+5,5) arc (0:350:5/18);
\draw[->] (18/18+5,5) arc (0:350:3/18);
\draw (10/9+5,5) -- (95/81+5,5);
\draw (100/81+5,5) -- (35/27+5,5);
\draw (40/27+5,5) arc (0:350:5/54);
\draw (40/27+5,5) -- (125/81+5,5);
\draw (130/81+5,5) -- (5/3+5,5);
\draw[->] (10/3+5,5) arc (0:350:5/6);
\draw[->] (19/6+5,5) arc (0:350:4/6);
\draw[->] (18/6+5,5) arc (0:350:3/6);
\draw (10/3+5,5) -- (275/81+5,5);
\draw (280/81+5,5) -- (95/27+5,5);
\draw (100/27+5,5) arc (0:350:5/54);
\draw (100/27+5,5) -- (305/81+5,5);
\draw (310/81+5,5) -- (35/9+5,5);
\draw[->] (40/9+5,5) arc (0:350:5/18);
\draw[->] (39/9+5,5) arc (0:350:3/18);
\draw (40/9+5,5) -- (365/81+5,5);
\draw (370/81+5,5) -- (125/27+5,5);
\draw (130/27+5,5) arc (0:350:5/54);
\draw (130/27+5,5) -- (395/81+5,5);
\draw (400/81+5,5) -- (5+5,5);
\end{tikzpicture}

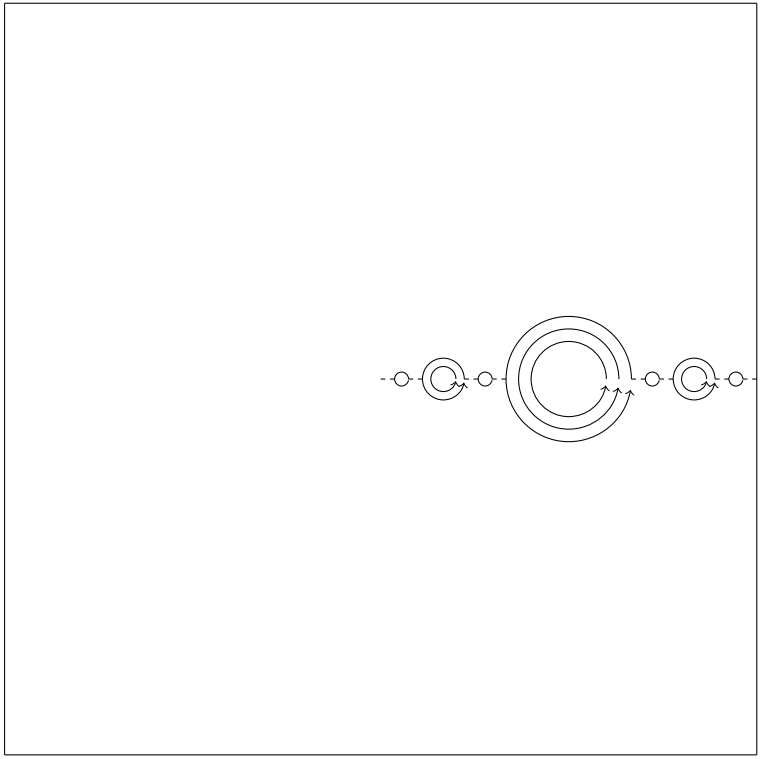
{\it Figure : Some illustration of the flow lines constructed in example 6.}

The vector field
$$   b(x,y)=\sum_{n=1}^\infty\sum_{ k=1,...,2^{n-1}}b_{n,k}(x,y) $$
generates circular flow lines around all  $ z_{n,k} $ for $ n\geq1 $ and $ 1\leq k\leq 2^{n-1} $.
It results from this that for any gap $ G_{n,k} $ its extremal points $ \inf G_{n,k} $ and $ \sup G_{n,k} $ are arbitrary close to some joint flow lines. which . Note that any two points $ u $ and $ v $ from $ C_n $ can be joint by a path which consists in using the flow to cross the gaps $ G_{n,k} $, which might separate $ u $ and $ v $. The maximal amount of cumulative jumps necessary to link $ u $ and $ v $ between those zero-cost flow portions is certainly less than the Lebesgue measure of $ C_n $. This implies that all points of the form $ (x,0);x\in C $ are in the same equivalence class, since the Lebesgue measure of $ C_n $ decreases to zero as $ n $ goes to infinity and $ C=\bigcap_{n\in\mathbb{N}}C_n $. 

Our goal is to define a continuous function $ g $ on the torus $ T^2 $ which is also in $ \overline{{\rm Ker}(B)\cap H^1} $ and which is constant on the flow-lines of the vector field $ b $ but not constant on the equivalence classes. In order to define the function $ g $, we put $ g(x,y)=f(x_{n,k}) $ for all $ (x,y)\in B_{n,k}:=B_{r_n}(z_{n,k}) $, where $ f $ is the Cantor function. As a result the function $ g $ is defined on the set 
$$ K=\bigcup_{n=1}^\infty\bigcup_{k=1}^{2^ {n-1}} B_{r_n}(z_{n,k}) .$$
Note that the function $ g $ is continuous on $ K $.

We now have to define the function $ g $ on the remaining parts of the torus. On the sets
$ \{(x,y);-1\leq x\leq 0\} $ we define $ g(x,y)=-x $. 

It then remains to define the function $ g $ on the set $ \{(x,y);0<x<1\}\backslash K $.  In the following we describe its construction on the upper part $ \{(x,y);0<x,y<1\}\backslash K $. The construction on the lower part $ \{(x,y);0<x<1, -1<y<0\}\backslash K $ carries over symmetrically. We first construct the familly of rectangles
\begin{eqnarray*}
  D_{1,1} &=& ]0,1[\times]r_1,1[ \\
  D_{2,1} &=& ]0,x_{1,1}-r_1[\times]r_2,r_1[ \\
  D_{2,2} &=& ]x_{1,1}+r_1,1[\times]r_2,r_1[ \\
  D_{3,1} &=& ]0,x_{2,1}-r_2[\times]r_3,r_2[ \\
  D_{3,2} &=& ]x_{2,1}+r_2,x_{1,1}-r_1[\times]r_3,r_2[ \\
  D_{3,3} &=& ]x_{1,1}+r_1,x_{2,2}-r_2[\times]r_3,r_2[ \\
 \rm{etc.} &&
\end{eqnarray*}
The following graphic illustrates the sets $ D_{n,k} $. Observe that for $ n\geq 1 $ one has
$$ \bigcup_{k=1}^{2^{n-1}}D_{n,k} = C_{n-1}\times]r_n,r_{n-1}[.
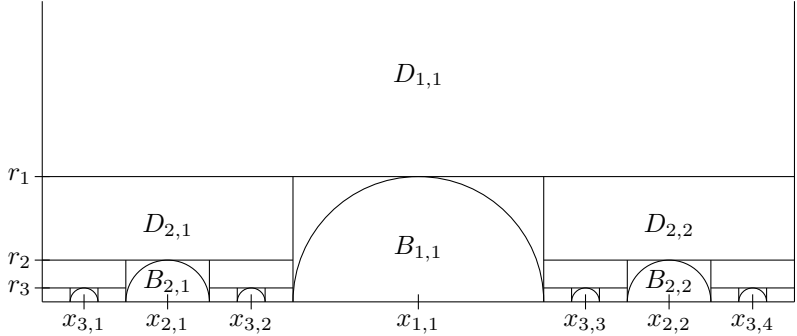
$$

\begin{tikzpicture}
\draw (0,0) -- (10,0);
\draw (0,0) -- (0,4);
\draw (10,0) -- (10,4);
\draw (40/6,0) arc (0:180:10/6);
\draw (40/18,0) arc (0:180:10/18);
\draw (160/18,0) arc (0:180:10/18);
\draw (40/54,0) arc (0:180:10/54);
\draw (160/54,0) arc (0:180:10/54);
\draw (400/54,0) arc (0:180:10/54);
\draw (520/54,0) arc (0:180:10/54);
\draw (0,10/6) -- (10,10/6);
\draw (40/6,0) -- (40/6,10/6);  
\draw (20/6,0) -- (20/6,10/6);  
\draw (20/18,0) -- (20/18,10/18);  
\draw (40/18,0) -- (40/18,10/18);  
\draw (140/18,0) -- (140/18,10/18);  
\draw (160/18,0) -- (160/18,10/18);  
\draw (20/54,0) -- (20/54,10/54);  
\draw (40/54,0) -- (40/54,10/54);  
\draw (140/54,0) -- (140/54,10/54);  
\draw (160/54,0) -- (160/54,10/54);  
\draw (380/54,0) -- (380/54,10/54);  
\draw (400/54,0) -- (400/54,10/54);  
\draw (500/54,0) -- (500/54,10/54);  
\draw (520/54,0) -- (520/54,10/54);  
\draw (0,10/18) -- (20/6,10/18);
\draw (40/6,10/18) -- (10,10/18);
\draw (0,10/54) -- (20/18,10/54);
\draw (40/18,10/54) -- (60/18,10/54);
\draw (120/18,10/54) -- (140/18,10/54);
\draw (160/18,10/54) -- (10,10/54);
\draw (10/2,-0.1) -- (10/2,0.1);
\draw (10/2,-0.3) node    {$x_{1,1}$};
\draw (10/2,+0.7) node    {$B_{1,1}$};
\draw (10/2,+3) node    {$D_{1,1}$};
\draw (10/6,-0.1) -- (10/6,0.1);
\draw (10/6,-0.3) node    {$x_{2,1}$};
\draw (10/6,+0.25) node    {$B_{2,1}$};
\draw (10/6,+1) node    {$D_{2,1}$};
\draw (50/6,-0.1) -- (50/6,0.1);
\draw (50/6,-0.3) node    {$x_{2,2}$};
\draw (50/6,+0.25) node    {$B_{2,2}$};
\draw (50/6,+1) node    {$D_{2,2}$};
\draw (10/18,-0.1) -- (10/18,0.1);
\draw (10/18,-0.3) node    {$x_{3,1}$};
\draw (50/18,-0.1) -- (50/18,0.1);
\draw (50/18,-0.3) node    {$x_{3,2}$};
\draw (130/18,-0.1) -- (130/18,0.1);
\draw (130/18,-0.3) node    {$x_{3,3}$};
\draw (170/18,-0.1) -- (170/18,0.1);
\draw (170/18,-0.3) node    {$x_{3,4}$};
\draw (-0.1,10/6) -- (0.1,10/6);
\draw (-0.3,10/6) node    {$ r_1 $};
\draw (-0.1,10/18) -- (0.1,10/18);
\draw (-0.3,10/18) node    {$ r_2 $};
\draw (-0.1,10/54) -- (0.1,10/54);
\draw (-0.3,10/54) node    {$ r_3 $};
\end{tikzpicture}

On the set $ D_{1,1} $ we define the function to satisfy $ f(x,y)=x $. 
It follows that $ g(x_{1,1}-r_1,r_1)=\tfrac{1}{2}-r_1 $ and $ g(x_{1,1}+r_1,r_1)=\tfrac{1}{2}+r_1 $. 
On the vertical line segment $ V_{1,1}^{(r)}=\big\{(x_{1,1}+r_1,y);0<y<r_1\big\} $ we then define $ g $ through linear interpolation as follows 
$$  g(x_{1,1}+r_1,y)=\tfrac{1}{2}+y .$$
Since by construction $ g(1,y)=1 $, we can define the values of $ g $ on $ D_{2,2} $ through horizontal linear interpolation as follows
$$   g(x,y)=(\tfrac{1}{2}+y)+\frac{x-\tfrac{1}{2}-r_1}{\tfrac{1}{2}-r_1}(\tfrac{1}{2}-y); \tfrac{1}{2}+r_1<x<1, r_2<y<r_1 . $$
A similar construction can be done to define $ g $ on the rectangle $ D_{2,1} $. 
As a result the partial derivatives of $ g $ on the rectangles $ D_{1,2} $ and $ D_{2,2} $ satisfy 
$$  \big|\partial_xg(x,y)\big|= \Big|\frac{\tfrac{1}{2}-y}{\tfrac{1}{2}-r_1}\Big| 
= \Big|\frac{\tfrac{1}{2}-y}{\tfrac{1}{2}-\tfrac{1}{2}\tfrac{1}{3}}\Big| \leq \Big|\frac{\tfrac{1}{2}}{\tfrac{1}{2}-\tfrac{1}{2}\tfrac{1}{3}}\Big| =\frac{3}{2} $$
and
$$ \big|\partial_yg(x,y)\big| = \Big|1-\frac{x-\tfrac{1}{2}-r_1}{\tfrac{1}{2}-r_1}\Big|
                                         = \Big|\frac{1-x}{\tfrac{1}{2}-r_1}\Big|
                                      \leq  \Big|\frac{1-\tfrac{1}{2}-r_1}{\tfrac{1}{2}-r_1}\Big| = 1 = (\tfrac{3}{2})^0 .$$
Now, one can iterate this process to define the value of the function $ g $ on all rectangles $ D_{n,k} $. For convenience we introduce a notation which helps to relate the set $ D_{n,k} $ to the center and radii of the balls which are to be found on its left and right side.
Note that for $ k>1 $ the rectangle $ D_{n,k} $ is delimited on its left side by a vertical line segment $ V_{n,k}^{(l)} $, which is tangential to one of the balls which forms the set $ K $.  We will refer to it as $ B_{n,k}^{(l)} $. Also le $ x_{n,k}^{(l)} $ be such that $ (x_{n,k}^{(l)},0) $ is the center of the ball $ B_{m,k}^{(l)} $. Further, denote by $ r_{n,k}^{(l)} $ the radius of the ball $ B_{n,k}^{(l)} $. For $ k=1 $ we put $ x_{n,1}^{(l)}=0 $ and $ r_{n,1}^{(l)}=0 $. The above line segment then can be represented as $ V_{n,k}^{(l)}=\big\{(x_{n,k}^ {(l)}+r_{n,k}^ {(l)},y);0<y<r_{n,k}^{(l)}\big\} $. In the same way, we use the fact that for $ k<2^{n-1} $ the set $ D_{n,k} $ is delimited on its right side by a vertical line segment $ V_{n,k}^{(r)} $ which is tangential to one of the balls $ B_{n,k}^{(r)} $ from $ K $ to define suitable $ x_{n,k}^{(r)} $ and $ r_{n,k}^{(r)} $. Furthermore, we put $ x_{n,2^{n-1}}^{(r)}=1 $, $ x_{n,1}^{(l)}=0 $ and $ r_{n,2^{n-1}}^{(r)}= r_{n,1}^{(l)}=0 $. The following figure illustrates the situation.

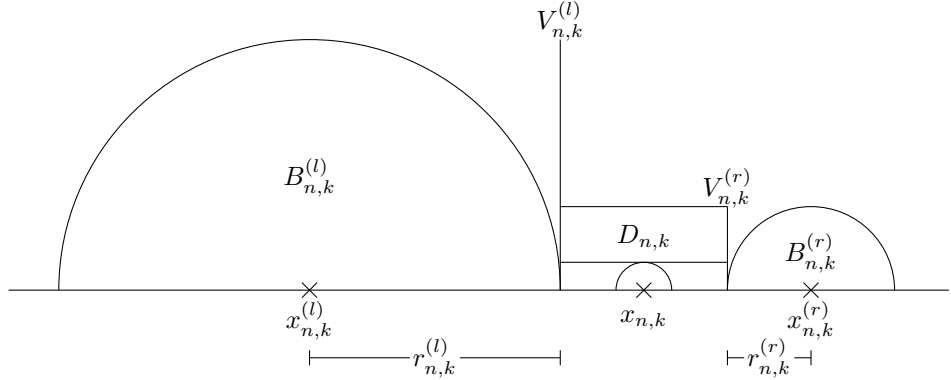
\begin{tikzpicture}
\draw (1,0) -- (13.5,0);
\draw (50/6,0) -- (50/6,20/6);
\draw (190/18,0) -- (190/18,20/18);
\draw (50/6,20/18) -- (190/18,20/18);
\draw (50/6,20/54) -- (190/18,20/54);
\draw (50/6,0) arc (0:180:20/6);
\draw (230/18,0) arc (0:180:20/18);
\draw (530/54,0) arc (0:180:20/54);
\draw (510/54,0.7) node    {$ D_{n,k} $};
\draw (210/18,0.5) node    {$ B_{n,k}^{(r)} $};
\draw (30/6,1.5) node    {$ B_{n,k}^{(l)} $};
\draw (210/18-0.1,-0.1) -- (210/18+0.1,0.1);
\draw (210/18-0.1,0.1) -- (210/18+0.1,-0.1);
\draw (30/6-0.1,-0.1) -- (30/6+0.1,0.1);
\draw (30/6-0.1,0.1) -- (30/6+0.1,-0.1);
\draw (210/18,-0.4) node    {$ x_{n,k}^{(r)} $};
\draw (30/6,-0.4) node    {$ x_{n,k}^{(l)} $};
\draw (510/54-0.1,-0.1) -- (510/54+0.1,0.1);
\draw (510/54-0.1,0.1) -- (510/54+0.1,-0.1);
\draw (510/54,-0.4) node  {$ x_{n,k} $};
\draw (50/6,20/6+0.25) node  {$ V_{n,k}^{(l)} $};
\draw (190/18,20/18+0.25) node  {$ V_{n,k}^{(r)} $};
\draw (50/6,-0.9) -- (42/6,-0.9);
\draw (38/6,-0.9) -- (30/6,-0.9);
\draw (50/6,-0.8) -- (50/6,-1);
\draw (30/6,-0.8) -- (30/6,-1);
\draw (40/6,-0.9) node  {$ r_{n,k}^{(l)} $};
\draw (190/18,-0.9) -- (194/18,-0.9);
\draw (206/18,-0.9) -- (210/18,-0.9);
\draw (190/18,-0.8) -- (190/18,-1);
\draw (210/18,-0.8) -- (210/18,-1);
\draw (200/18,-0.9) node  {$ r_{n,k}^{(r)} $};
\end{tikzpicture}

Now, in the process of defining the function $ g $ on the set $ D_{n,k} $ we already know its values on the horizontal line segment $ H_{n,k}=\big]x_{n,k}^{(l)}+r_{n,k}^{(l)},x_{n,k}^{(r)}-r_{n,k}^{(r)}\big[\times \{r_{n-1}\} $. 
Thus its values in the upper endpoints of the vertical line segments $ V_{n,k}^{(l)} $ and $ V_{n,k}^{(r)} $ are known. 
The values of $ g $ on those line segments are then obtained through linear interpolation between those values and $ f_{n,k}^{(l)}:=f(x_{n,k}^{(l)}) $ resp. $ f_{n,k}^{(r)}:=f(x_{n,k}^{(r)}) $. As a consequence on $ V_{n,k}^{(l)} $ the function $ g $ has the following representation
$$ g(x_{n,k}^{(l)}+r_{n,k}^{(l)},y)=f_{n,k}^{(l)}+\delta_{n,k}^{(l)}y  \ \ \ \mbox{for $ 0<y<r_{n,k}^{(l)} $}  $$ 
and on $ V_{n,k}^{(l)} $ the function $ g $ has the expression
$$ g(x_{n,k}^{(r)}-r_{n,k}^{(r)},y)= f_{n,k}^{(r)}-\delta_{n,k}^{(r)}y \ \ \ \mbox{for $ 0<y<r_{n,k}^{(r)} $}   ,$$
with suitable positive coefficients $ \delta_{n,k}^{(l)} $ and $ \delta_{n,k}^{(r)} $. Those coefficients correspond to the rate of increase of the function $ x\mapsto g(x,r_{n,k}^{(l)}) $ resp. $ x\mapsto g(x,r_{n,k}^{(r)}) $. As we will see below, they are upper bounded by $ (\tfrac{3}{2})^{n-2} $.
The function $ g $ is then defined on the rectangle $ D_{n,k} $ through horizontal interpolation as follows
$$   g(x,y)=\Big((f_{n,k}^{(r)}-\delta_{n,k}^{(r)}y)-(f_{n,k}^{(l)}+\delta_{n,k}^{(l)}y)\Big)\frac{x-x_{n,k}^{(l)}-r_{n,k}^{(l)}}{x_{n,k}^{(r)}-r_{n,k}^{(r)}-x_{n,k}^{(l)}-r_{n,k}^{(l)}} + (f_{n,k}^{(l)}+\delta_{n,k}^{(l)}y) .$$
Note that by construction of the Cantor function, one has $ f_{n,k}^{(r)}-f_{n,k}^{(l)}=2^{n-1} $. Moreover, by construction of $ D_{n,k} $ one has $ x_{n,k}^{(r)}-r_{n,k}^{(r)}-x_{n,k}^{(l)}-r_{n,k}^{(l)}=3^{n-1} $. It thus follows 
\begin{eqnarray*}
 \big|\partial_xg(x,y)\big| 
 &=& \Bigg|\frac{(f_{n,k}^{(r)}-\delta_{n,k}^{(r)}y)-(f_{n,k}^{(l)}+\delta_{n,k}^{(l)}y)}{x_{n,k}^{(r)}-r_{n,k}^{(r)}-x_{n,k}^{(l)}-r_{n,k}^{(l)}}\Bigg| \\
&=& 3^{n-1}\Big|(f_{n,k}^{(r)}-\delta_{n,k}^{(r)}y)-(f_{n,k}^{(l)}+\delta_{n,k}^{(l)}y)\Big|\\
&\leq& 3^{n-1}\Big|f_{n,k}^{(r)}-f_{n,k}^{(l)}\Big| \\
 &\leq& (\tfrac{3}{2})^{n-1}.
\end{eqnarray*}
By the symmetry of the construction the rate of increase of the function $ g $ along the segments $ V_{n,k}^{(l)} $ reproduces its rate of increase along the horizontal line segment $ \big]x_{n,k}^{(l)}-r_{n,k}^{(l)},x_{n,k}^{(l)}+r_{n,k}^{(l)}\big[\times\big\{r_{n,k}^{(l)}\big\} $. Thus it follows that
$$  0<\delta_{n,k}{(l)} \leq \big|\partial_xg(x,r_{n,k}^{(l)})\big| \leq (\tfrac{3}{2})^{m-1} \ \ \ \mbox{for $ 1\leq m<n $.} $$
Thus it follows that $ 0<\delta_{n,k}^{(l)}<(\tfrac{3}{2})^{n-2} $. Since a similar reasoning also applies to $ \delta_{n,k}^{(r)} $, we have for all $ 1\leq k\leq 2^{n-1} $ that
$$  0<\min(\delta_{n,k}^{(l)},\delta_{n,k}^{(r)})<\max(\delta_{n,k}^{(l)},\delta_{n,k}^{(r)})\leq (\tfrac{3}{2})^{n-2} .$$
Further, one has
\begin{eqnarray*}
 \big|\partial_yg(x,y)\big| 
 &=&\Bigg|\big(\delta_{n,k}^{(r)}-\delta_{n,k}^{(l)}\big)\frac{x-x_{n,k}^{(l)}-r_{n,k}^{(l)}}{x_{n,k}^{(r)}-r_{n,k}^{(r)}-x_{n,k}^{(l)}-r_{n,k}^{(l)}}+\delta_{n,k}^{(l)}\Bigg| \\
&\leq&\big|\delta_{n,k}^{(r)}-\delta_{n,k}^{(l)}\big|\Bigg|\frac{x-x_{n,k}^{(l)}-r_{n,k}^{(l)}}{x_{n,k}^{(r)}-r_{n,k}^{(r)}-x_{n,k}^{(l)}-r_{n,k}^{(l)}}\Bigg| +\delta_{n,k}^{(l)}\\
&\leq& 2\max\big(\delta_{n,k}^{(r)},\delta_{n,k}^{(l)}\big) \\
 &\leq&  2(\tfrac{3}{2})^{n-2} .
\end{eqnarray*}
We so far have seen that on $ D_{n,k} $ we have
$$ \big|\nabla g\big|^2=(\partial_xg)^2+(\partial_yg)^2\leq   \Big((\tfrac{3}{2})^{n-1}\Big)^2+\Big(2(\tfrac{3}{2})^{n-2}\Big)^2=3(\tfrac{3}{2})^{2n-2} .$$
It thus follows that
\begin{eqnarray*}
 \int_{\bigcup_{n,k}D_{n,k}} \big|\nabla g(x,y)\big|^2dxdy 
 &=& 3\sum_{n=1}^\infty\sum_{k=1}^{2^{n-1}}\int_{D_{n,k}}(\tfrac{3}{2})^{2n-2}dxdy \\
 &\leq& 3\sum_{n=1}^\infty\sum_{k=1}^{2^{n-1}}(\tfrac{3}{2})^{2(n-1)}\lambda^{\otimes 2}(D_{n,k}) \\
 &=& 3\sum_{n=1}^\infty\sum_{k=1}^{2^{n-1}}(\tfrac{3}{2})^{2(n-1)}\big(r_{n-1}-r_n\big)(\tfrac{1}{3})^{n-1} \\
 &=& 3\sum_{n=1}^\infty2^{n-1}(\tfrac{3}{2})^{2(n-1)}\Big(\tfrac{1}{2}(\tfrac{1}{3})^{n-1}-\tfrac{1}{2}(\tfrac{1}{3})^n\Big)(\tfrac{1}{3})^{n-1} \\
&=& \sum_{n=1}^\infty2^{n-1}(\tfrac{3}{2})^{2(n-1)}(\tfrac{1}{3})^{n-1}(\tfrac{1}{3})^{n-1} \\
&=& \sum_{n=1}^\infty(\tfrac{1}{2})^{n-1} \ < \ \infty.
\end{eqnarray*}

We now turn our attention to the remaining part of $ ]0,1[^2\backslash K $. For this, we define the sets
$$  R_{n,k}=\Big\{(x,y): x_{n,k}<x<x_{n,k}+r_{n}, \sqrt{r_{n}^2-x^2}<y<r_{n}\Big\} $$
and
$$  L_{n,k}=\Big\{(x,y): x_{n,k}-r_{n}<x<x_{n}, \sqrt{r_{n}^2-x^2}<y<r_{n}\Big\} $$
Those two sets can also be represented in polar coordinates around the point $ (x_{n,k},1) $ as follows
$$  R_{n,k}=\Big\{(\varphi,r); 0<\varphi<\tfrac{\pi}{2}, r_{n}<r<r_{n}/\cos\varphi\Big\} $$

The remaining part, were the function $ g $ has to be defined is the union of disjoint sets
$$  \bigcup_{n=1}^\infty\bigcup_{k=1}^{2^{n-1}}L_{n,k}\cup R_{n,k} . $$
In order to define the function $ g $ on those two sets, we use polar coordinates around the point $ (x_{n,k},1) $ and introduce the functions $ r(\varphi) $ and $ h(\varphi) $ according to the follwing figure.\\

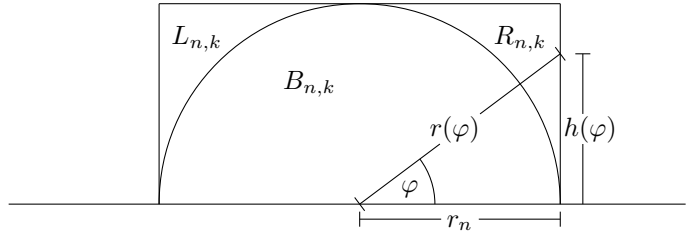
\begin{tikzpicture}
\draw (1,0) -- (10,0);
\draw (50/6,0) -- (50/6,16/6);
\draw (18/6,0) -- (18/6,16/6);
\draw (18/6,16/6) -- (50/6,16/6);
\draw (34/6,0) -- (34/6+16/6*2/5,0+12/6*2/5);
\draw (34/6+16/6*3/5,0+12/6*3/5) -- (34/6+16/6,0+12/6);
\draw (50/6-2/30,2+8/90) -- (50/6+2/30,2-8/90);
\draw (34/6-2/30,8/90) -- (34/6+2/30,-8/90);
\draw (50/6,0) arc (0:180:16/6);
\draw (40/6,0) arc (0:37:1);
\draw (38/6,0.2) node    {$ \varphi $};
\draw (30/6,1.6) node    {$ B_{n,k} $};
\draw (42/6,-0.25) node    {$ r_{n} $};
\draw (3.5,2.2) node    {$ L_{n,k} $};
\draw (7.8,2.2) node    {$ R_{n,k} $};
\draw (34/6+16/6*2/5+0.2,0+12/6*2/5+0.2) node    {$ r(\varphi) $};
\draw (34/6,-0.2) -- (40.5/6,-0.2);
\draw (43.5/6,-0.2) -- (50/6,-0.2);
\draw (50/6,-0.1) -- (50/6,-0.3);
\draw (34/6,-0.1) -- (34/6,-0.3);
\draw (50/6+0.3,0) -- (50/6+0.3,0.8);
\draw (50/6+0.3,1.2) -- (50/6+0.3,2);
\draw (50/6+0.2,2) -- (50/6+0.4,2);
\draw (50/6+0.4,1) node    {$ h(\varphi) $};
\end{tikzpicture}

We then have 
$$   r(\varphi)=\frac{r_n}{\cos\varphi} \ \ \ \mbox{and} \ \ \ h(\varphi)= r(\varphi)\sin\varphi=r_n\tan\varphi .$$
Note that the values of $ g $ are already defined on the boundary of the region $ R_{n,k} $.
Its values inside of the set $ R_{n,k} $ are then obtained through interpolation along lines passing through $ x_{n,k} $:
The value on the set $ B_{n,k} $ equals $ f_{n,k}:=f(x_{n,k}) $. We already discussed that on the vertical segment $ \{x_{n,k}+r_{n}\}\times]0,r_{n}[ $ which forms the right boundary of $ R_{n,k} $ we have $ g(x_{n,k}+r_{n},y)=\delta_{n,k}y+f_{n,k} $ with some suitable $ \delta_{n,k} $ satisfying $ |\delta_{n,k}|\leq(\frac{3}{2})^{n-1} $. Further, on the horizontal segment $ ]x_{n,k},x_{n,k}+r_{n}[\times\{r_{n}\} $ forming the bottom boundary of $ R_{n,k} $ we have $ g(x,y)=\delta_{n,k}(x+x_{n,k})+f_{n,k} $.
By symmetry it is sufficient to study the interpolation of the function along straight lines intersecting the vertical segment $ \{x_{n,k}+r_{n}\}\times]0,r_{n}[ $.
We there then define the function $ g $ with respect to polar coordinates as
\begin{eqnarray*}
  g(r,\varphi) &=& f_{n,k}+\frac{r_{n}-r}{r_n-r(\varphi)}(g(x_{n,k}+r_{n},h(\varphi)) -f_{n,k})  \\
  &=& f_{n,k}+\frac{r_n-r}{r_{n}-\frac{r_n}{\cos\varphi}}(g(x_{n,k}+r_{n},r_n\tan\varphi)-f_{n,k})   \\
  &=& f_{n,k}+\frac{r_n-r}{r_{n}-\frac{r_n}{\cos\varphi}}\delta_{n,k}r_n\tan\varphi   \\
  &=& f_{n,k}+(r_n-r)\frac{\sin\varphi}{\cos\varphi-1}\delta_{n,k}  .
\end{eqnarray*}
Its partial derivative with respect to $ r $ equals
\begin{eqnarray*}
  \partial_r g(r,\varphi)  &=& -\delta_{n,k}\frac{\sin\varphi}{\cos\varphi-1}
\end{eqnarray*}
and its partial derivative with respect to $ \varphi $ equals
\begin{eqnarray*}
  \partial_\varphi g(r,\varphi)  &=& \delta_{n,k}(r_n-r)\frac{\cos\varphi(\cos\varphi-1)+\sin\varphi\sin\varphi}{(\cos\varphi-1)^2} \\
     &=&   \delta_{n,k}(r_n-r)\frac{1-\cos\varphi}{(\cos\varphi-1)^2}\\
   &=&  \delta_{n,k}(r_n-r)\frac{1}{1-\cos\varphi}.
\end{eqnarray*}
It follows that
\begin{eqnarray*}
 \iint_{R_{n,k}}|\partial_r g|^2dxdy 
   &=&\int_0^{\frac{\pi}{4}}\int_{r_{n}}^{r(\varphi)}r|\partial_r g(r,\varphi)|^2drd\varphi \\
    &=& \delta_{n,k}^2\int_0^{\frac{\pi}{4}}\int_{r_{n}}^{r(\varphi)}r\frac{\sin^2\varphi}{(\cos\varphi-1)^2}drd\varphi \\
    &=& \delta_{n,k}^2\int_0^{\frac{\pi}{4}}\tfrac{1}{2}(r(\varphi)^2-r_{n}^2)\frac{\sin^2\varphi}{(\cos\varphi-1)^2}d\varphi \\
  &\leq& \delta_{n,k}^2\tfrac{1}{2}r_{n}^2(2-1)\int_0^{\frac{\pi}{4}}\frac{\sin^2\varphi}{(\cos\varphi-1)^2}d\varphi  \ < \ K_1\delta_{n,k}^2r_{n}^2
\end{eqnarray*}
and also since $ r(\varphi)-r_{n}=\frac{r_n}{\cos\varphi}-r_n=r_n(\frac{1}{\cos\varphi}-1)=r_n\frac{1-\cos\varphi}{\cos\varphi} $ 
\begin{eqnarray*}
 \iint_{R_{n,k}}|\tfrac{1}{r}\partial_\varphi g|^2dxdy 
    &=&\int_0^{\frac{\pi}{4}}\int_{r_{n}}^{r(\varphi)}r|\tfrac{1}{r}\partial_\varphi g(r,\varphi)|^2drd\varphi \\
&=&\delta_{n,k}^2\int_0^{\frac{\pi}{4}}\int_{r_{n}}^{r(\varphi)}\frac{1}{r}\frac{(r-r_{n})^2}{(\cos\varphi-1)^2}drd\varphi \\
&\leq&\delta_{n,k}^2\frac{1}{r_{n}}\int_0^{\frac{\pi}{4}}\int_{r_{n}}^{r(\varphi)}\frac{(r(\varphi)-r_{n,k})^2}{(\cos\varphi-1)^2}drd\varphi \\
 &=& \delta_{n,k}^2\frac{1}{r_{n}}\int_0^{\frac{\pi}{4}}\frac{(r(\varphi)-r_{n})^3}{(\cos\varphi-1)^2}d\varphi \\
 &=& \delta_{n,k}^2r_{n}^2\int_0^{\frac{\pi}{4}}\frac{(1-\cos\varphi)}{\cos^3\varphi}d\varphi 
  \ < \ K_2\delta_{n,k}^2r_{n}^2 .
\end{eqnarray*}
As a consequence we have that
$$  \iint_{R_{n,k}}|\nabla g|^2dxdy 
    \ = \ \iint_{R_{n,k}}|\partial_r g|^2dxdy +\iint_{R_{n,k}}|\tfrac{1}{r}\partial_\varphi g|^2dxdy 
     \ \leq \ K\delta_{n,k}^2r_{n}^2 . $$
It follows from this that
\begin{eqnarray*}
 \iint_{\bigcup_{n,k}R_{n,k}}|\nabla g|^2dxdy 
  &=& \sum_{n=1}^\infty\sum_{k=1}^{2^{n-1}}\iint_{R_{n,k}}|\nabla g|^2dxdy \\
  &\leq& K\sum_{n=1}^\infty\sum_{k=1}^{2^{n-1}}\delta_{n,k}^2r_{n}^2 \\
  &\leq& K\sum_{n=1}^\infty\sum_{k=1}^{2^{n-1}}\big((\tfrac{3}{2})^{n-1}\big)^2\big(\tfrac{1}{2}(\tfrac{1}{3})^n\big)^2 \\
  &\leq& K(\tfrac{1}{6})^2\sum_{n=1}^\infty 2^{n-1}(\tfrac{3}{2})^{2n-2}(\tfrac{1}{3})^{2n-2}\\
  &=& K(\tfrac{1}{6})^2\sum_{n=1}^\infty 2^{n-1}(\tfrac{1}{2})^{2n-2} \ < \ \infty .
\end{eqnarray*}
All in all we have shown that
$$  \iint_{\mathbb{T}^2}|\nabla g|^2dxdy <\infty .$$
Thus the function $ g $ is in $  H^1\cap C(\mathbb{T}^2) $ and constant along all trajectories. However, it is not constant on the equivalence classes defined in section \ref{Sektion5}.

\section{Appendix}

{\bf Proof of Proposition \ref{Rellich} :} 
The first part is a special case of Riesz theorem on the weak compactness of the unit-ball in a separable Hilbert-space. 
We give a proof for the convenience of the reader. 
Assume that there exists a sequence $ \psi_n,n\geq1 $ in $ H^1 $ such that $ \|\psi_n\|_{H^1}^2<C $ for all $ n\in\N $.
There exist suitable coefficients $ a_{n,k},n\in\N,k\in\N $ such that 
$$  \psi_n=\sum_{k=1}^\infty a_{n,k}\phi_k\ \ \ 
   {\rm and}\ \ \ \sum_{k=1}^\infty a_{n,k}^2\lambda_k<C .$$
This implies that the sequence $ a_{n,1},n\in\N $ satisfies $ a_{n,1}^2\lambda_1<C $. 
There exists an infinite set $ \tilde{N}_1\subset \N $ and an $ a_1\in\R $ such that 
$$  a_{n,1}\rightarrow a_1 \ {\rm as}\ n\rightarrow\infty\ {\rm in}\ \tilde{N}_1.$$
We note that we can find an infinite subset $ N_1\subset\tilde{N}_1 $ such that  
$$ \sum_{k=2}^\infty a_{n,k}^2\lambda_k\leq C-a_1^2\lambda_1+1/2 \ \ \ {\rm for\ all}\ n\in N_1.$$
We now observe that the sequence $ a_{n,2},n\in N_1 $ also satisfies $ a_{n,2}^2\lambda_2<C $.
If we iterate the above procedure, we find a sequence of infinite subsets $ N_m;m\in\N $ with $ N_m\subset N_{m-1} $ and a sequence of real numbers 
$ a_m;m\in\N $ such that 
$$   a_{n,m}\rightarrow a_m \ {\rm as}\ n\rightarrow\infty\ {\rm in}\ N_m $$ 
and 
$$  \sum_{k=m}^\infty a_{n,k}^2\lambda_k\leq C-\sum_{l=1}^{m-1}a_l^2\lambda_l+(1/2)^m
                    \ \ \ {\rm for\ all}\ n\in N_m .$$
If we choose a diagonal-sequence with the property $ n(m)\in N_m $ for all $ m\in\N $, then it follows that for all $ k\in\N $ we have
$$ a_{n(m),k}\rightarrow a_k \ \ \ {\rm as}\ m\rightarrow\infty $$
and 
$$  \sum_{k=m}^\infty a_{n(m),k}^2\lambda_k\leq C-\sum_{l=1}^{m-1}a_l^2\lambda_l+(1/2)^m
                    \ \ \ {\rm for\ all}\ m\in\mathbb{N} .$$
This last statement implies that $ \psi :=\sum_{k=1}^\infty a_k\phi_k $ is in $ H^1 $ since we also have
$$  \sum_{k=1}^\infty a_k^2\lambda_k\leq C .$$                              
Now, we fix an arbitrary $ \varphi\in H^1 $. It then can be represented in the following form
$$ \varphi=\sum_{k=1}^\infty c_k\phi_k \ \ \ 
   {\rm with}\ \ \ \sum_{k=1}^\infty c_k^2\lambda_k<\infty .$$
For an arbitrary $ \epsilon>0 $ there exists a $ K\in\N $ such that 
$$    \sum_{k>K}c_k^2\lambda_k<\epsilon^2/(C+1) .$$
Moreover, we have for all $ m\in\N $ that 
\begin{eqnarray*}
 \sum_{k>K}|a_{n(m),k}||c_k|\lambda_k&\leq& \left(\sum_{k>K}a_{n(m),k}^2\lambda_k\right)^{1/2}\left(\sum_{k>K}c_k^2\lambda_k\right)^{1/2}<\epsilon .
\end{eqnarray*}
In the same way we have 
\begin{eqnarray*}
  \sum_{k>K}|a_k||c_k|\lambda_k&\leq& \left(\sum_{k>K}a_k^2\lambda_k\right)^{1/2}\left(\sum_{k>K}c_k^2\lambda_k\right)^{1/2}<\epsilon .
\end{eqnarray*}
This yields that
\begin{eqnarray*}
 &&\left| \langle\psi_{n(m)},\varphi\rangle_{H^1}-\langle\psi,\varphi\rangle_{H^1}\right|\\
 &\leq&\sum_{k>K}|a_{n(m),k}||c_k|\lambda_k
   +\sum_{k\leq K}|a_{n(m),k}-a_k||c_k|\lambda_k+\sum_{k>K}|a_k||c_k|\lambda_k\\
 &\leq&2\epsilon+\sum_{k\leq K}|a_{n(m),k}-a_k||c_k|\lambda_k.
\end{eqnarray*}
However, the last term converges to zero as $ m\rightarrow\infty $, because 
$ a_{n(m),k}\rightarrow a_k $ for all $ k\in\N $. This proves the first part, since $ \epsilon $ can be 
arbitrarily small. \\[1mm]
The second part is a type of Sobolev-embedding theorem for our situation, since it states that the inclusion
of $ H^1 $ into $ H $ is compact. Here we use the fact that $ \lambda_n\rightarrow\infty $.
We have to prove that $ \psi_{n(m)} $ converges toward $ \psi $ in $ H $. 
For all $ N\in\mathbb{N} $ and $ m\in\mathbb{N} $  we have 
\begin{eqnarray*}
   \sum_{k=N}^\infty |a_{n(m),k}-a_k|^2&\leq&\frac{1}{\lambda_N}\sum_{k>N}|a_{n(m),k}-a_k|^2\lambda_k\\
   &\leq& \frac{1}{\lambda_N}\|\psi_{n(m)}-\psi\|_{H^1}^2
   \leq\frac{4C}{\lambda_N}.
\end{eqnarray*}
Thus there exists for all $ \epsilon>0 $ a $ N>0 $ such that 
$$  \sum_{k>N}^\infty |a_{n(m),k}-a_k|^2<\epsilon .$$
This implies that for all $ m\in\mathbb{N} $
\begin{eqnarray*}
  \|\psi_{n(m)}-\psi\|^2=\sum_{k=1}^\infty |a_{n(m),k}-a_k|^2\leq\epsilon+\sum_{k\leq N}|a_{n(m),k}-a_k|^2.
\end{eqnarray*}
The last term converges toward zero as $ m\rightarrow\infty $. This proves the result, since $ \epsilon $ can 
be chosen arbitrarily small. 
\hfill $ \Box $

\section*{Acknowledgements}
This research was made possible through several visits of Brice Franke at National Central University and Academia Sinica in Taipei where supported by NCTS and the Research Grant ...


\begin{thebibliography}{xxx}


\bibitem{BvR}{\scshape F.~Barret}, {\scshape M.~von Renesse}  : Averaging principle for diffusion processes via Dirichlet forms, 
{\slshape Potential Analysis}, {\bfseries 41}, 1033-1063, (2014).

\bibitem{BLP}{\scshape A.~Bensoussan}, {\scshape J.L.~Lions}, {\scshape G.~Papanicolaou}  : {\slshape Asymptotic Analysis for Periodic Structures}, North-Holland, (1978).

\bibitem{BGW}{\scshape R.~Bhattacharya}, {\scshape V.~Gupta}, {\scshape H.~Walker}  : Asymptotics of solute dispersion in periodic porous media, 
{\slshape SIAM Journal on Applied Mathematics}, {\bfseries 49}, 86-98, (1989).

\bibitem{Bil} {\scshape P.~ Billingsley}  : {\slshape Convergence of Probability Measures}, John Wiley \& Sonc, Inc. New York, (1968).

\bibitem{Cha} {\scshape I.~Chavel}  : {\slshape Eigenvalues in Riemannian Geometry}, Academic Press Inc., Orlando, (1984).

\bibitem{CKRZ} {\scshape P.~Constantin}, {\scshape A.~Kiselev}, {\scshape L.~Ryzhik}, {\scshape A.~Zlatos}  : Diffusion and mixing in fluid flow, 
{\slshape Annals of Mathematics}, {\bfseries 168}, 643-674, (2008).

\bibitem{DFK}  {\scshape D.~Dolgopyat},  {\scshape M.~Freidlin},  {\scshape L.~Koralov}  : Deterministic and stochastic perturbations of area preserving flows on a two-dimensional torus,
{\slshape Ergodic Theory and Dynamical Systems}, {\bfseries 32}, 899-918, (2012).

\bibitem{Dieud} {\scshape J.~Dieudonn\'{e}}  : {\slshape Foundations of modern analysis}, Vol.1, Academic Press, New York-London, (1969).

\bibitem{DLP}{\scshape A. Duncan, T. Leli\`evre and G. Pavliotis }  : Variance reduction using nonreversible langevin samplers.   {\slshape Journal of Statistical Physics}, 163, 457-49, (2016).

\bibitem{FHPS} {\scshape B.~Franke}, {\scshape C.-R.~Hwang}, {\scshape H.-M.~Pai}, {\scshape S.-J.~Sheu}  : The behavior of the spectral gap under growing drift, {\slshape Transactions of the  American Mathematical Society}, {\bfseries 362}, 1325-1350, (2010).

\bibitem{FHOP}{\scshape B.~Franke}, {\scshape C.-R.~Hwang}, , {\scshape A.~Ouled-Said}, {\scshape H.-M.~Pai} : A note on the asymptotic variance of drift accelerated diffusion, {\slshape Statistics and Probability Letters}, {\bfseries 175}, 109128, (2021).

\bibitem{FW} {\scshape M.~Freidlin}, {\scshape M.~Weber}  : Random perturbations of dynamical systems and diffusion processes with conservation
laws, {\slshape Probability Theory and Related Fields}, {\bfseries 128}, 441-466, (2004).

\bibitem{Fukushima} {\scshape M.~Fukushima}  : {\slshape Dirichlet Forms and Markov Processes}, North Holland Publishing Company, Amsterdam, (1980).

\bibitem{IsSo} {\scshape H.~Ishii}, {\scshape P.~Souganidis}  : A PDE approach to small stochastic perturbations of Hamiltonian flows, 
{\slshape Journal of Differential Equations}, {\bfseries 252}, 1748-1775,  (2012).

\bibitem{Helf} {\scshape B.~Helffer}  : {\slshape Semiclassical Analysis, Witten Laplacians and Statistical Mechanics}, World Scientific, (2002).

\bibitem{Hsu} {\scshape E.~Hsu}  : {\slshape Stochastic Analysis on Manifolds}, GSM 38, AMS, (2001).

\bibitem{HNW}{\scshape C.-R. Hwang,  R. Normand, and S.-J. Wu }  :  Variance reduction for diffusions.  {\slshape Stochastic Processes and their Applications}, 125(9), 3522-3540, (2015).

\bibitem{HS} {\scshape C.-R.~Hwang}, {\scshape S.-J.~Sheu}  : On some quadratic perturbation of Ornstein-Uhlenbeck processes, {\slshape Soochow Journal of Mathematics}, {\bfseries 26}, 205-244, (2000).

\bibitem{Kling} {\scshape W.~Klingenberg}  : {\slshape Riemannian Geometry}, de Gruyter, (1995).

\bibitem{KarShr} {\scshape I.~Kartzas}, {\scshape S.~Shreve}  : {\slshape Brownian Motion and Stochastic Calculus}, GTM 113, Springer Verlag New-York, (1987).

\bibitem{Kre} {\scshape U. Krengel}  : {\slshape Ergodic Theorems}, Studies in Mathematics, 6, Walter de Gruyter, Berlin, (1985).

\bibitem{Man} {\scshape R.~Man\'{e}}  : {\slshape Ergodic Theory and Differentiable Dynamics}, Springer Verlag Berlin, (1987).

\bibitem{Paz} {\scshape A.~Pazy}  : {\slshape Semigroups of Linear Operators and Applications to Partial Differential Equations},
Applied Mathematical Sciences, 44, Springer Verlag New-York, (1983).

\bibitem{Yoshida} {\scshape K.~Yoshida}  : {\slshape Functional Analysis}, Springer Verlag Berlin, (1980).


\end{thebibliography}
\end{document}